\documentclass[11 pt]{amsart}

\usepackage{amssymb} \usepackage{amsfonts} \usepackage{amsmath}
\usepackage{amsthm} \usepackage{epsfig} 

\newtheorem{lemma}{Lemma}[section]
\newtheorem{teo}[lemma]{Theorem}
\newtheorem{prop}[lemma]{Proposition}
\newtheorem{cor}[lemma]{Corollary}

\theoremstyle{definition}
\newtheorem{defn}[lemma]{Definition}

\newtheorem{quest}[lemma]{Question}

\newtheorem{example}[lemma]{Example}

\theoremstyle{remark}
\newtheorem{rem}[lemma]{Remark}

\newcommand{\matR} {\ensuremath {\mathbb{R}}}
\newcommand{\matQ} {\ensuremath {\mathbb{Q}}}
\newcommand{\matZ} {\ensuremath {\mathbb{Z}}}

\newcommand{\cerchio}{\includegraphics[width = .4 cm]{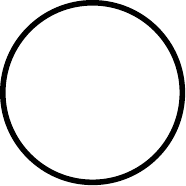}}
\newcommand{\teta}{\includegraphics[width = .4 cm]{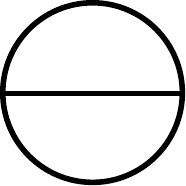}}
\newcommand{\tetra}{\includegraphics[width = .4 cm]{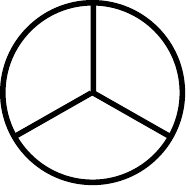}}
\newcommand{\tetracolored}{\includegraphics[width = 1 cm] {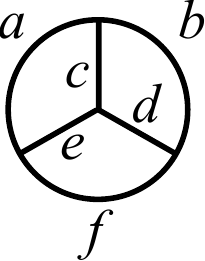}}
\newcommand{\tetracoloredtwo}{\includegraphics[width = 1 cm]{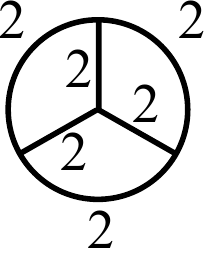}}

\newcommand{\gl}{{\rm gl}}
\newcommand{\ord}{{\rm ord}}

\newcommand{\dimo}[1]{\vspace{2pt}\noindent\textit{Proof of \ref{#1}}.\ }
\newcommand{\finedimo}{{\hfill\hbox{$\square$}\vspace{2pt}}}

\author{Alessio Carrega}
\address{Dipartimento di Matematica ``Tonelli'', Largo Pontecorvo 5, 56127 Pisa, Italy}
\email{carrega at mail dot dm dot unipi dot it}

\author{Bruno Martelli}
\address{Dipartimento di Matematica ``Tonelli'', Largo Pontecorvo 5, 56127 Pisa, Italy}
\email{martelli at dm dot unipi dot it}

\title{Shadows, ribbon surfaces, \\ and quantum invariants}

\begin{document}

\begin{abstract}
Eisermann has shown that the Jones polynomial of a $n$-component ribbon link $L\subset S^3$ is divided by the Jones polynomial of the trivial $n$-component link. 
We improve this theorem by extending its range of application from links in $S^3$ to colored knotted trivalent graphs in $\#_g(S^2\times S^1)$, the connected sum of $g\geqslant 0$ copies of $S^2\times S^1$. 

We show in particular that if the Kauffman bracket of a knot in $\#_g(S^2\times S^1)$ has a pole in $q=i$ of order $n$, the ribbon genus of the knot is at least $\frac {n+1}2$. We construct some families of knots in $\#_g(S^2\times S^1)$ for which this lower bound is sharp and arbitrarily big. We prove these estimates using Turaev shadows.


\end{abstract}

\maketitle

\setcounter{tocdepth}{1}
\tableofcontents

\section{Introduction} \label{introduction:section}
Thirty years after its discovery, we know only a few relations between the Jones polynomial $J_L$ of a link $L$ and its topological properties. A notable one is Eisermann's Theorem \cite{Eis} which connects the Jones polynomial to four-dimensional smooth topology. The theorem states that the Jones polynomial of a $n$-component ribbon link is divided by the Jones polynomial of the trivial $n$-component link. 

Another four-dimensional object related to the Jones polynomial is Turaev's \emph{shadow}. In this paper we reprove Eisermann's Theorem using shadows, and extend its range of application from links in $S^3$ to colored trivalent graphs in $\#_g(S^2\times S^1)$. 

In this introduction, we first show how we re-prove Eisermann's theorem for links in $S^3$, and later explain its extension to graphs in $\#_g(S^2\times S^1)$.


\subsection{Shadows}
Shadows are simple two-dimensional polyhedra locally-flatly embedded in four-manifolds. They were defined by Turaev \cite{Tu:preprint, Tu} and then considered by various authors, see for instance \cite{Bu, C:complexity, CoThu, CoThu:preprint, Go, IK, M, Shu, Thu, Tu:paper}. 

In this paper, a shadow $X$ is a (simple, locally-flat) \emph{collapsible spine} of $D^4$. 
Being a collapsible spine is a quite restrictive requirement: we want 
\begin{enumerate}
\item that $D^4$ collapses to $X$ (\emph{i.e.}~$X$ is a spine of $D^4$), 
\item that $X$ collapses to a point (\emph{i.e.}~$X$ is collapsible).
\end{enumerate}
If we use the symbols $\searrow$ and $\bullet$ to indicate collapsing and a point, we may summarize that by writing
$$D^4 \ \searrow \ X\ \searrow\ \bullet$$

Recall that a \emph{ribbon surface} is a properly embedded surface $S\subset D^4$ that can be put into Morse position with only minima and saddle points (no maxima). The surface $S$ may be disconnected and non-orientable. We will start by proving the following purely topological fact: 

\begin{teo} \label{main:topological:teo}
Every ribbon surface is contained in some shadow.
\end{teo}

\begin{figure}
\begin{center}
\includegraphics[width = 8 cm]{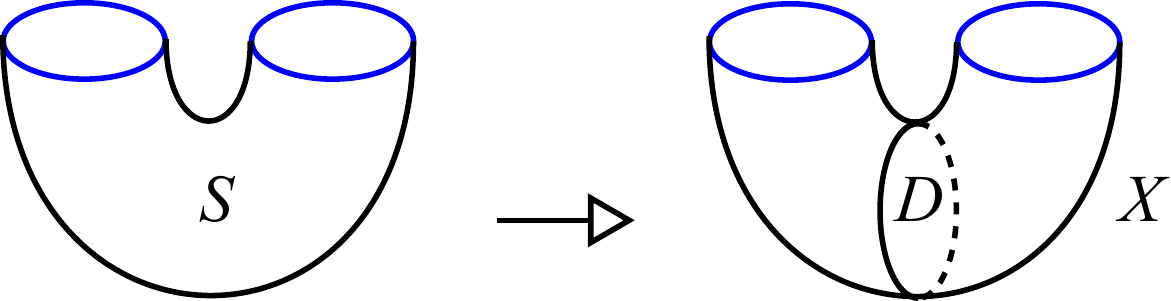}
\caption{If $S\subset D^4$ is the trivially embedded annulus bounding the unlink $\partial S$, a shadow $X= S \cup D$ is obtained by attaching a disc $D$ to its core.}
\label{shadow_annulus_0:fig}
\end{center}
\end{figure}

We single out a couple of examples.

\begin{example}
Consider the trivially embedded annulus $S\subset D^4$ as in Fig.~\ref{shadow_annulus_0:fig}. A shadow $X$ containing $S$ is constructed by attaching a disc $D$ to its core. Note that indeed $D^4 \searrow X \searrow \bullet$. See Example \ref{annulus:example}.
\end{example}

\begin{example} The trivial ribbon disc (with one minimum and no saddles) is itself a shadow of $D^4$. However, a non-trivial ribbon disc $D\subset D^4$ is \emph{not} a shadow: it collapses to a point, but it fails to be a spine of $D^4$, see Proposition \ref{trivial:prop}. The shadow containing $D$ may be rather complicate.
\end{example}

The hypothesis that the surface is ribbon is crucial here: there are surfaces (for instance, discs) that are not contained in any shadow. Indeed the following implications hold for a properly embedded surface $S\subset D^4$:
\begin{align*} 
S {\rm \ ribbon\ } \Longrightarrow \ S{\rm\ contained} {\rm\ in\ a\ shadow}\ \Longrightarrow\ S {\rm\ homotopically\ ribbon}.
\end{align*}
Recall that $S$ is \emph{homotopically ribbon} if the inclusion 
$S^3\setminus \partial S \hookrightarrow D^4\setminus S$
induces a surjective homomorphism on fundamental groups. It is easy to construct discs that are not homotopically ribbon, and hence are not contained in any shadow, see Section \ref{non-ribbon:subsection}.

The question whether every homotopically ribbon surface $S$ is actually ribbon is, to the best of our knowledge, open: the requirement that $S$ is contained in some shadow lies between these two properties and breaks this question in two parts.

\begin{quest} \label{intro:quest}
Can we reverse any of the two implications above?
\end{quest}

\subsection{Quantum invariants}
We then turn to quantum invariants. A shadow $X$ for a link $L\subset S^3$ is a shadow $X\subset D^4$ such that $X\cap S^3 = L$. An easy homological argument shows that $X$ contains a unique surface $S$ with $\partial S = L$. The surface $S$ is possibly disconnected and non-orientable, but it contains no closed components.

Instead of the Jones polynomial $J_L$ we prefer to use the Kauffman bracket $\langle L \rangle$ that is more adapted to our purposes. The Kauffman bracket is a Laurent polynomial in $q$ that differs from $J_L$ only by some re-parametrization.

The shadow $X$ can be used to calculate $\langle L\rangle$ via Turaev's state-sum formula \cite{Tu:preprint, Tu}, and by analyzing carefully that formula we prove the following:

\begin{teo} \label{main:quantum:teo}
Let $X$ be a shadow for a link $L\subset S^3$ and $S\subset X$ be the unique surface with $\partial S = L$. The Kauffman bracket $\langle L \rangle $ vanishes at least $\chi(S)$ times at $q=i$.
\end{teo}

The theorem provides some information only when $\chi(S) > 0$.
The two theorems we stated imply Eisermann's Theorem \cite{Eis}:

\begin{cor}
If a link $L\subset S^3$ bounds a ribbon surface $S$ then $\langle L\rangle$ vanishes at least $\chi(S)$ times at $q=i$.
\end{cor}
\begin{proof}
There is a shadow $X$ for $L$ that contains $S$ by Theorem \ref{main:topological:teo}. Theorem \ref{main:quantum:teo} implies that $\langle L\rangle$ vanishes at least $\chi(S)$ times at $q=i$.
\end{proof}

Again, this corollary is relevant only when $\chi(S)>0$. Recall that a $n$-component link $L\subset S^3$ is \emph{ribbon} if it bounds a ribbon surface that consists of $n$ discs. The interesting corollary is of course the following.

\begin{cor}
If $L\subset S^3$ is a $n$-component ribbon link then $\langle L \rangle$ vanishes at least $n$ times at $q=i$.
\end{cor}

It is an immediate consequence of its definition that the Kauffman bracket $\langle L \rangle$ of any link $L\subset S^3$ vanishes at $q=i$ at least once, and hence Eisermann's Theorem actually provides some information only when $n\geqslant 2$. In particular, unfortunately it says nothing about knots in $S^3$. 

\subsection{Links in some other manifolds}
The techniques we used to prove Theorems \ref{main:topological:teo} and \ref{main:quantum:teo} extend naturally in two directions.
The first consists in varying the ambient manifold.

Costantino \cite{C:Jones} has defined the Kauffman bracket $\langle L \rangle$ of a framed link 
$$L\subset \#_g(S^2\times S^1)$$ 
in a connected sum of any $g\geqslant 1$ copies of $S^2\times S^1$. The Kauffman bracket $\langle L \rangle$ is now a rational function on $q^{\frac 12}$ which may have poles at some roots of unity, including the value $q=i$ we are interested in. So we define 
$$\ord_i\langle L \rangle \in \matZ \cup \{+\infty\}$$ 
to be
the maximum integer $k$ such that $\frac{\langle L \rangle}{(q-i)^{k-1}}$ vanishes in $q=i$.  This is the first exponent of the Laurent expansion of $\langle L \rangle$ at $q=i$.

The notion of ribbon surface extends naturally to any closed 3-manifold $M$: a ribbon surface is a properly embedded surface in Morse position inside $M\times [0,1]$, with boundary in $M\times 0$, and without maxima. Equivalently, it is an immersed surface in $M$ having only ribbon singularities, see Section \ref{ribbon:subsection}. We generalize Eisermann's theorem as follows:

\begin{teo} \label{link:teo}
If a link $L\subset \#_g(S^2\times S^1)$ bounds a ribbon surface $S$ then 
$$\ord_i\langle L \rangle \geqslant \chi(S).$$
\end{teo}

This theorem is potentially stronger in $\#_g(S^2\times S^1)$ than in $S^3$ because now $\ord_i\langle L \rangle$ can be an arbitrarily small \emph{negative} number. In particular it provides non-trivial informations also for knots, as the following example shows.

\begin{figure}
\begin{center}
\includegraphics[width = 10 cm]{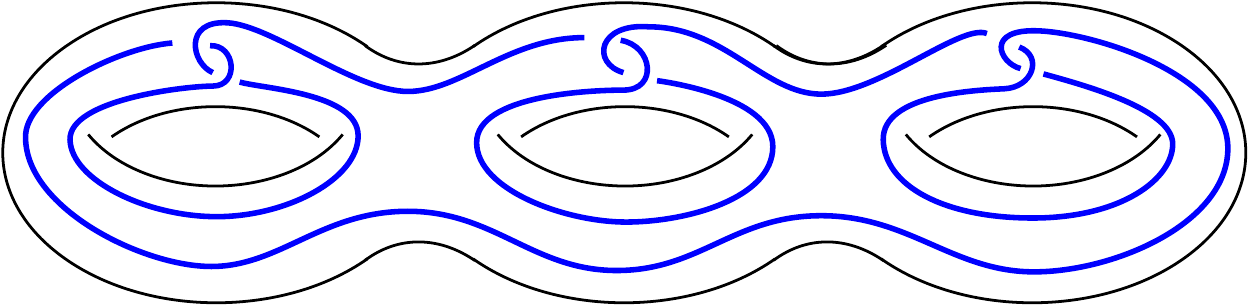}
\caption{A knot $K$ in $\#_g(S^2\times S^1)$. To get $\#_g(S^2\times S^1)$ simply double the handlebody in the picture. We draw here the case $g=3$, the general case is obvious from the picture. Note that the knot is null-homotopic.}
\label{knot_ribbon:fig}
\end{center}
\end{figure}

\begin{example} \label{nodo:example}
The framed knot $K\subset \#_g(S^2\times S^1)$ drawn in Fig.~\ref{knot_ribbon:fig} has 
$$\langle K \rangle =  (-1)^{1-g} q^{-\frac{3g}2}\frac{(1+ q^2+q^4+q^6)^g}{(q+q^{-1})^{2g-1}}$$
and hence $\ord_i\langle K \rangle = g-(2g-1) = 1-g$. Therefore $K$ bounds no ribbon surface $S$ with $\chi(S) >1-g$. In particular, it is not a ribbon knot.
\end{example}

The \emph{ribbon genus} of a knot $K$ is the minimum genus of an orientable connected ribbon surface $S$ with $\partial S = K$. As a consequence, the ribbon genus of the knot $K$ shown in Fig.~\ref{knot_ribbon:fig} is at least $\frac {g}2$. In general we get:

\begin{cor}
Let $K\subset \#_g(S^2\times S^1)$ be a knot. Then:
\begin{itemize}
\item if $\langle K \rangle$ does not vanish at $q=i$, the knot is not ribbon,
\item if $\langle K \rangle$ has a pole at $q=i$ of order $n>0$, the ribbon genus of $K$ is at least $\frac {n+1}2$.
\end{itemize}
\end{cor}

The following example illustrates a family of knots for which this lower bound on the ribbon genus is sharp and arbitrarily big.

\begin{example} \label{easy:example}
Let $S$ be a compact orientable surface with boundary. The boundary of the four-manifold $S\times D^2$ is diffeomorphic to $\#_g(S^2\times S^1)$ for $g=1-\chi(S)$. The link $L=\partial S$ inside $\#_g(S^2\times S^1)$ has
$$\langle L \rangle = (-q-q^{-1})^{\chi(S)}$$
and hence $\ord_i\langle L \rangle = \chi(S)$.
Therefore $S$ is a ribbon surface of maximal Euler characteristic (among those having $L$ as boundary). The lower bound given by Theorem \ref{link:teo} is sharp on these links. We can choose $L$ to be a knot by picking a surface $S$ with one boundary component, and we can choose $\chi(S)$ to be arbitrarily small by increasing the genus of $S$.
\end{example}

\begin{rem}
Similar lower bounds for the slice genus of the knots and links considered in Examples \ref{nodo:example} and \ref{easy:example} can be constructed by other methods, see Remark \ref{easy:was:rem}: these basic examples were chosen primarily because their Kauffman bracket can be easily calculated by hand. 

Since the lower bound furnished by the Kauffman bracket is non-trivial on these simple examples, it might hopefully say something relevant on more elaborate ones: we briefly discuss the slice/ribbon conjecture and its possible extensions in Section \ref{slice-ribbon:subsection}.
\end{rem}

\subsection{Knotted trivalent graphs}
The second extension consists of taking trivalent graphs instead of just links.
The Kauffman bracket $\langle G \rangle$ is defined for colored framed knotted trivalent graphs $G$ in $S^3$ and more generally in $\#_g(S^2\times S^1)$. These objects are often called \emph{ribbon graphs}, but we do not use this terminology here to avoid confusion with ribbon surfaces.

The coloring of $G$ is the assignment of a non-negative integer to every edge or knot component of $G$, such that at every vertex $v\in G$ the colors $a,b,c$ of the incident edges fulfill the triangle inequalities and have even sum $a+b+c$. Thanks to these admissibility conditions, the numbers
$$\frac{a+b-c}2, \quad \frac{b+c-a}2, \quad \frac{c+a-b}2$$
are non-negative integers. We say that the vertex $v$ is \emph{red} if at least two of these integers are odd. The edges in $G$ having an odd color form a sublink $L\subset G$ called the \emph{odd sublink}.

The Kauffman bracket $\langle G \rangle$ of $G$ is still a rational function in $q^{\frac 12}$. The following theorem generalizes Theorem \ref{link:teo} from links to graphs.

\begin{teo} \label{large:teo}
Let $G$ be a colored framed knotted trivalent graph in $S^3$ or $\#_g(S^2\times S^1)$ and $L\subset G$ be its odd sublink. If $L$ bounds a ribbon surface $S$ then
$$\ord_i\langle G \rangle \geqslant \chi(S) - \frac r2$$
where $r$ is the number of red vertices in $G$.
\end{teo}

The theorem applies in particular to colored links:

\begin{cor} Let $G$ be a colored framed link in $S^3$ or $\#_g(S^2\times S^1)$ and $L\subset G$ be its odd sublink. If $L$ bounds a ribbon surface $S$ then
$$\ord_i\langle G \rangle \geqslant \chi(S).$$
\end{cor}

Hence in particular Eisermann's Theorem holds \emph{as is} for links colored with odd numbers.

\subsection{Proofs}
The proof of Theorem \ref{large:teo} splits into two parts: the topological Theorem \ref{main:topological:teo}, and the more technical Theorem \ref{main:quantum:teo}, both extended from links in $S^3$ to graphs in $\#_g(S^2\times S^1)$. 

While the topological side of the story is a one-page proof, the technical part needs a long case-by-case analysis that we would have never pursued if we were not aware of Eisermann's Theorem. We easily \emph{localize} the proof of Theorem \ref{large:teo} to the case where $G$ is one of the three planar graphs 
$$\cerchio, \teta, \tetra$$ 
in $S^3$. The graph $\tetra$ is a well-known building block in quantum topology (closely related to the \emph{quantum $6j$-symbols}) and its Kauffman bracket is a quite complicate rational function in $q$, see Section \ref{basic:subsection}.

To prove Theorem \ref{large:teo} we examine carefully this rational function near $q=i$ for all possible parities of the six numbers coloring the edges of the graph, and check that the inequality $\ord_i \tetra \geqslant \chi(S) - \frac{r}2$ is fulfilled (quite miraculously) in all cases (and it is almost always an equality!). The addendum $\frac r2$ in the formula is absolutely necessary, as the following shows.

\begin{example}
The Kauffman bracket of the graph $G = \teta_{2,2,2}$ colored with $2,2,2$ is
$$\langle G \rangle = -\frac{(q^3+q+q^{-1}+q^{-3})(q^2+1+q^{-2})}{(q+q^{-1})^2}$$
and has a pole in $q=i$ of order 1, \emph{i.e.}~$\ord_i(G)=-1$. The odd sublink of $G$ is empty and hence bounds the empty ribbon surface $S$ that has $\chi(S)=0$. The formula $\ord_i(G) \geqslant \chi(S) - \frac r2$ holds because both vertices of $G$ are red and hence $r=1$, giving $-1 \geqslant 0-1$.
\end{example}

\subsection{Structure of the paper}
We define ribbon surfaces and shadows in Section \ref{topological:section}, where we also prove the topological Theorem \ref{main:topological:teo}. In Section \ref{Jones:section} we introduce the Kauffman bracket and recall Turaev's formula that computes it as a state-sum on a shadow. In Section \ref{estimate:section} we prove the more technical Theorem \ref{main:quantum:teo}. In Section \ref{other:section} we generalize everything from $S^3$ to $\#_g(S^2\times S^1)$.
In Section \ref{state-sum:section} we re-prove Turaev's state-sum formula.
Section \ref{questions:section} is devoted to some open questions for further research.

\subsection{Acknowledgements}
We would like thank Francesco Costantino, Paolo Lisca, and Dylan Thurston for many helpful conversations.

\section{Shadows and ribbon surfaces} \label{topological:section}

We introduce ribbon surfaces and shadows, and then prove Theorem \ref{main:topological:teo} which says that every ribbon surface is contained in some shadow.

\subsection{Ribbon surfaces} \label{ribbon:subsection}
A properly embedded smooth surface $S\subset D^4$ is \emph{ribbon} if one of the following equivalent conditions holds:

\begin{figure}
\begin{center}
\includegraphics[width = 4 cm]{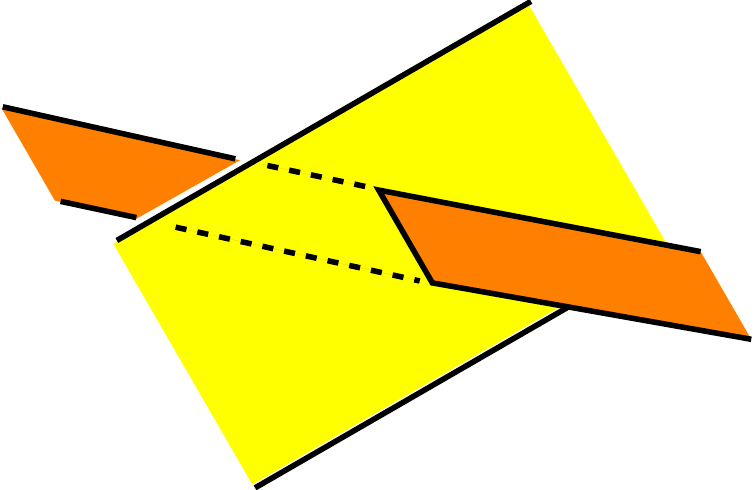}
\caption{A ribbon singularity}
\label{ribbon_singularity:fig}
\end{center}
\end{figure}

\begin{figure}
\begin{center}
\includegraphics[width = 9 cm]{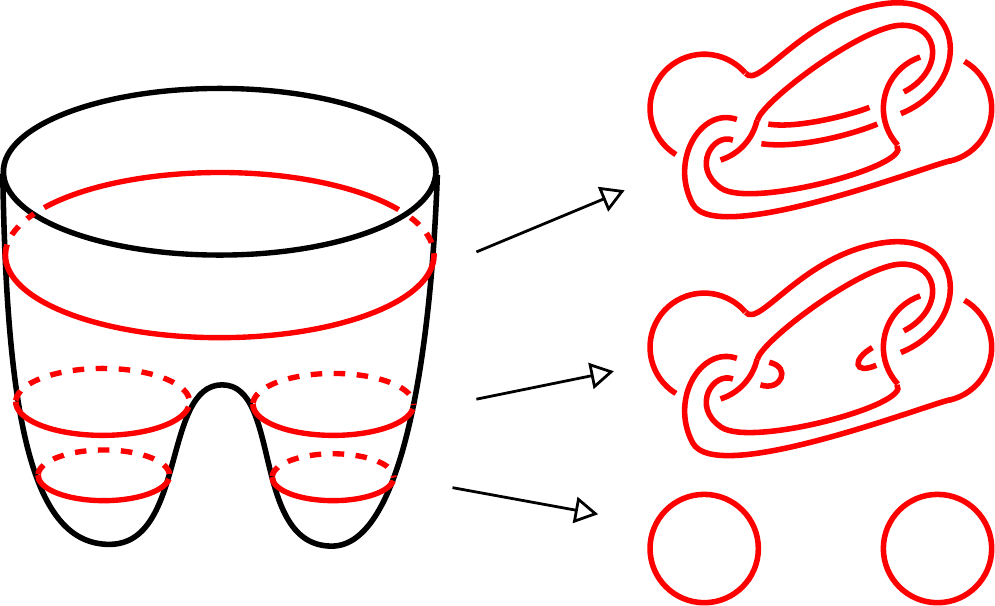}
\caption{A ribbon disc in $D^4$ in Morse position with two minima and one saddle. Each regular level gives a link in $S^3$.}
\label{Morse:fig}
\end{center}
\end{figure}

\begin{itemize}
\item the surface $S$ may be isotoped to an immersed surface in $S^3$ having only ``ribbon'' singularities as in Fig.~\ref{ribbon_singularity:fig},
\item the surface $S$ may be isotoped in $D^4$ into Morse position, with only minima and saddle points (no maxima) as in Fig.~\ref{Morse:fig}.
\end{itemize}

\begin{figure}
\begin{center}
\includegraphics[width = 12.5 cm]{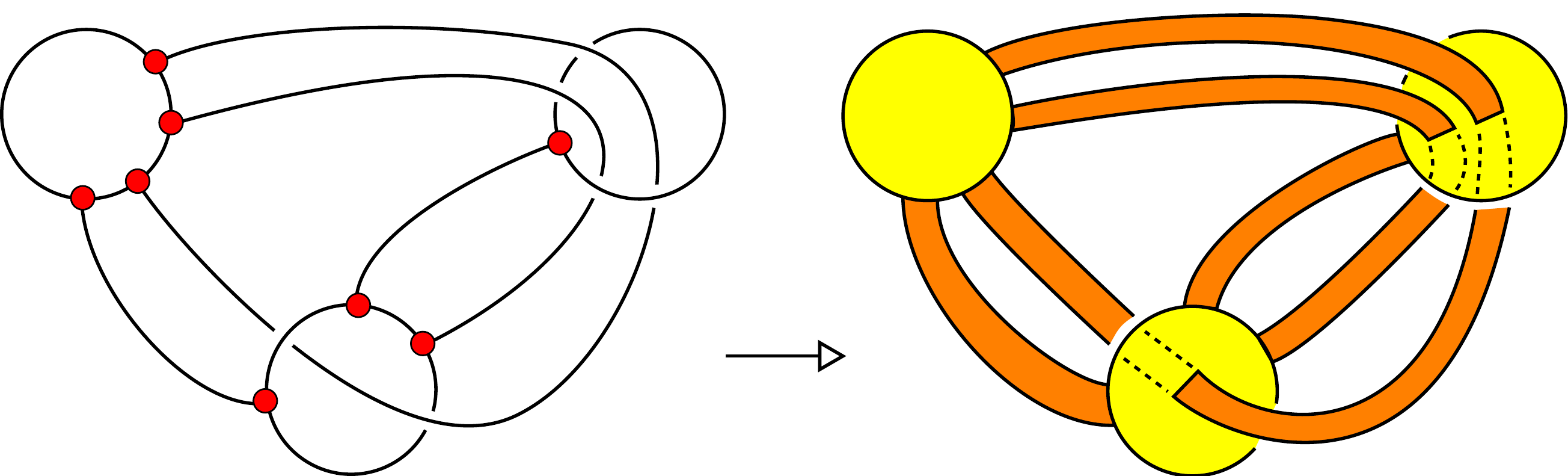}
\caption{Every ribbon surface can be constructed from a planar diagram with some disjoint circles representing the minima and some edges connecting them representing the saddles (left). The surface is obtained by filling the circles (yellow) and thickening the edges to (orange) bands.}
\label{construct_ribbon:fig}
\end{center}
\end{figure}

Every ribbon surface $S$ can be constructed from a planar diagram as in Fig.~\ref{construct_ribbon:fig}, consisting of some disjoint circles and some arcs connecting them in space.

\subsection{Shadows}

\begin{figure}
\begin{center}
\includegraphics[width = 12 cm]{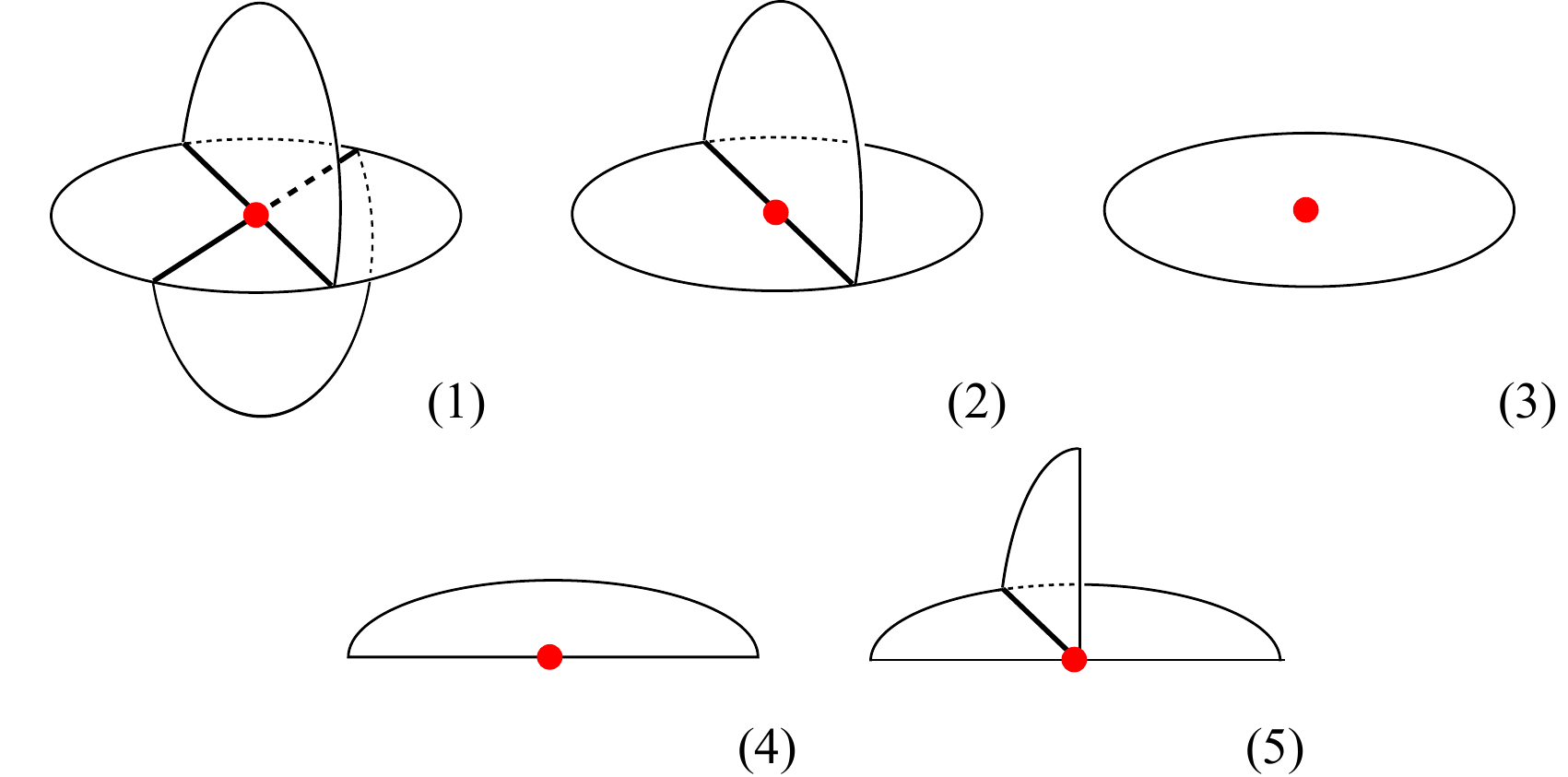}
\caption{Neighborhoods of points in a simple polyhedron.}
\label{models:fig}
\end{center}
\end{figure}

A \emph{simple polyhedron} $X$ is a 2-dimensional compact polyhedron where every point has a neighborhood homeomorphic to one the five types (1-5) shown in Fig.~\ref{models:fig}. The five types form subsets of $X$ whose connected components are called \emph{vertices} (1), \emph{interior edges} (2), \emph{regions} (3), \emph{boundary edges} (4), and \emph{boundary vertices} (5). The points (4) and (5) altogether form the \emph{boundary} $\partial X$ of $X$. An edge is either an open segment or a circle; a region is a (possibly non-orientable) connected surface.

\begin{defn}
A \emph{shadow} for $D^4$ is a simple polyhedron $X\subset D^4$ such that the following holds:
\begin{itemize}
\item $X$ is properly embedded, that is $\partial X = X\cap S^3$,
\item $X$ is locally flat: every point $p\in X$ has a neighborhood $U$ in $D^4$ diffeomorphic to $B^3 \times (-1,1)$ with $U\cap X$ contained in $B^3 \times 0$ as in Fig.~\ref{models:fig},
\item $X$ collapses onto a point,
\item $D^4$ collapses onto $X$.
\end{itemize}
\end{defn}

The first two conditions are just reasonable requirements one assumes when considering simple polyhedra inside four-manifolds; on the other hand, the last two conditions are quite restrictive and can be summarized by writing
$$D^4 \ \searrow \ X\ \searrow\ \bullet$$
where $\bullet$ indicates a point. 

\subsection{Knotted trivalent graphs}
A \emph{knotted trivalent graph} is a smooth graph in $S^3$ where every vertex has valence $3$, and knot components are also admitted. So in particular a link is a knotted trivalent graph without vertices. 

The boundary $G=\partial X$ of a shadow $X\subset D^4$ is a knotted trivalent graph in $S^3$, and we say that $X$ is a shadow of $G$. Although the definition of shadow seems very restrictive, it turns out that every knotted trivalent graph has at least one shadow (and in fact, infinitely many):

\begin{figure}
\begin{center}
\includegraphics[width = 11 cm]{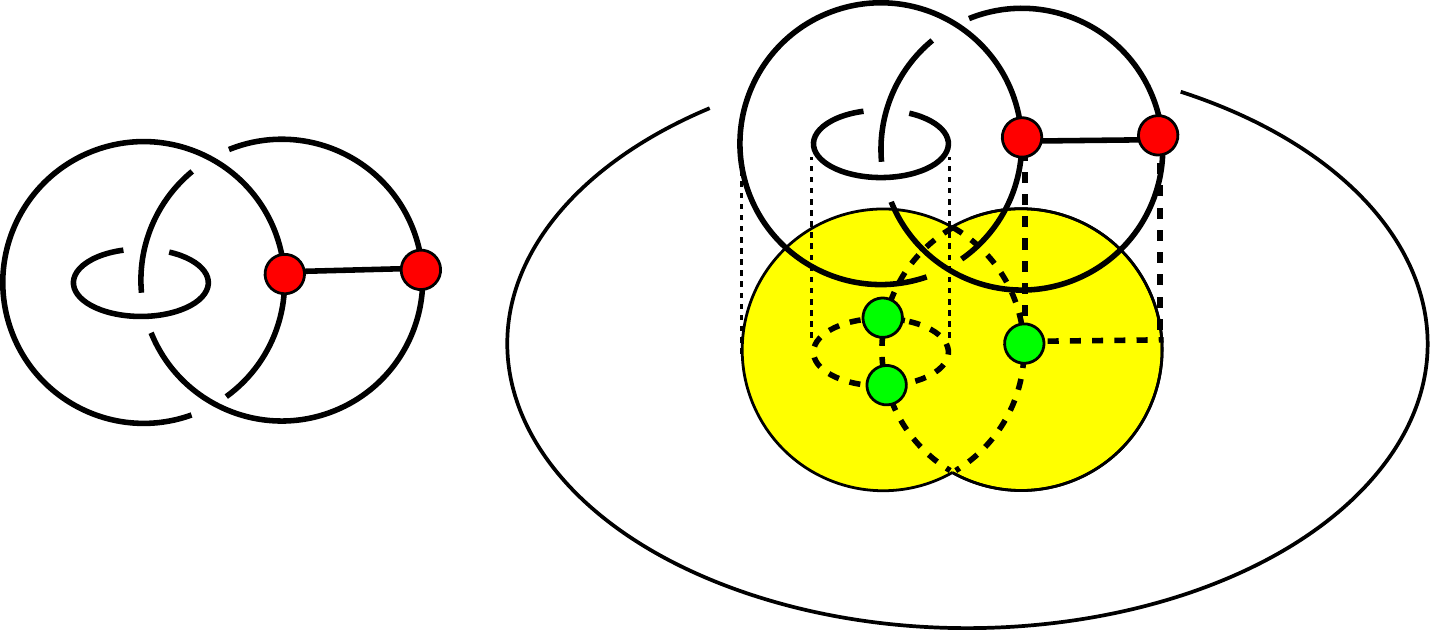}
\caption{A knotted trivalent graph (left) and its shadow (right). The shadow contains three interior vertices (green).}
\label{KTG:fig}
\end{center}
\end{figure}

\begin{prop}[Turaev] \label{Turaev:prop}
Every knotted trivalent graph $G\subset S^3$ has a shadow.
\end{prop}
\begin{proof}
This result was first proved by Turaev \cite{Tu:preprint} in a more general context; here we follow the proof contained in \cite[Theorem 3.14]{CoThu}. Pick a diagram for $G$ as in Fig.~\ref{KTG:fig}-(left). We suppose that there is a smallest closed disc $D$ containing the diagram like the yellow one in Fig.~\ref{KTG:fig}-(right). This is equivalent to ask that the diagram is connected and no vertex of the diagram disconnects it: these conditions can be easily achieved using Reidemeister moves.

If we push the yellow disc $D$ entirely inside $D^4$, we obviously get $D^4 \searrow D \searrow \bullet$. We enlarge $D$ by adding a cylinder above $G$ as sketched in Fig.~\ref{KTG:fig}-(right). The resulting object $X$ is a shadow for $G$: we still have $D^4 \searrow X \searrow \bullet$, and $X$ is easily seen to be a properly embedded locally flat simple polyhedron with $\partial X = G$.
\end{proof}

\subsection{Ribbon surfaces in shadows}
We are ready to prove Theorem \ref{main:topological:teo}:

\begin{teo} \label{ribbon:contained:teo}
Every ribbon surface $S$ is contained in a shadow $X$ with $\partial X = \partial S$.
\end{teo}
\begin{proof}
Construct $S$ from a planar diagram $G$ as in Fig.~\ref{construct_ribbon:fig}-(left). Via Reidemeister moves we may suppose that there is a smallest closed disc containing $G$. The diagram $G$ identifies a knotted trivalent graph and we construct a shadow $X$ for $G$ using the algorithm described in the proof of Proposition \ref{Turaev:prop}.

Note that $X$ contains the yellow discs of Fig.~\ref{construct_ribbon:fig}-(right). To complete the construction, we simply add to $X$ the orange bands shown in Fig.~\ref{construct_ribbon:fig}-(right), and then push their interior a bit inside $D^4$. We end up with a shadow $X$ containing the whole of $S$ and with $\partial X = \partial S$.
\end{proof}

\begin{figure}
\begin{center}
\includegraphics[width = 11 cm]{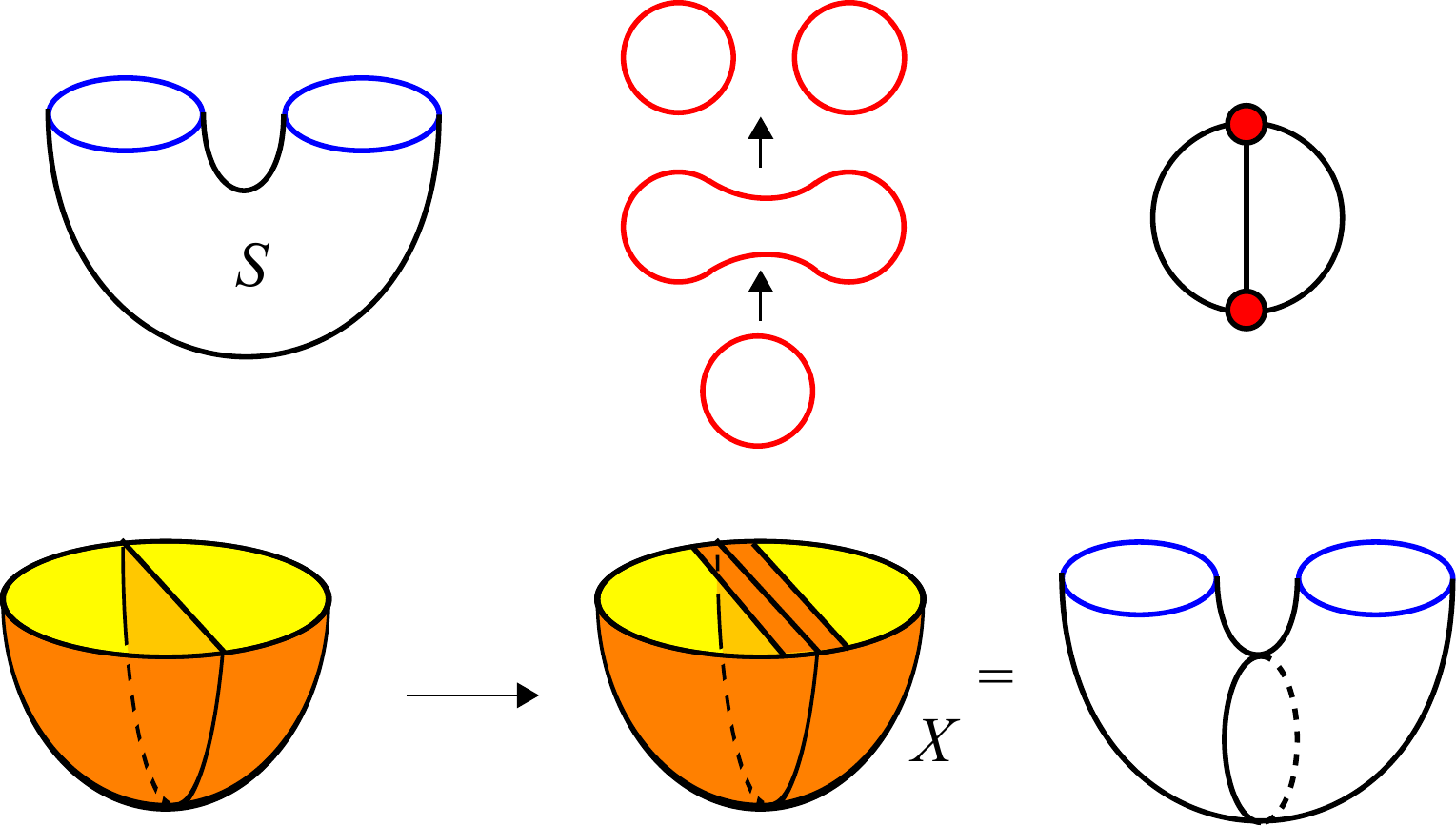}
\caption{How to build a shadow $X$ containing a given ribbon surface $S$. We show here the construction for the ribbon annulus $S$.}
\label{shadow_annulus:fig}
\end{center}
\end{figure}

\begin{example} \label{annulus:example}
Fig.~\ref{shadow_annulus:fig} illustrates the construction in a simple case. The ribbon surface $S\subset D^4$ is a trivially embedded annulus with boundary $L=\partial S$ the unlink with two components; the annulus $S$ in Morse position has one minimum and one saddle, and it is hence a ribbon surface constructed from the graph $G$ shown in Fig.~\ref{shadow_annulus:fig}-(top-right): a circle (the minimum) with a diameter (encoding the saddle). A shadow for $G$ is shown in Fig.~\ref{shadow_annulus:fig}-(bottom-left). By adding a band we obtain a shadow $X$ for $L$ containing $S$, and $X$ is just $S$ with a disc attached to its core. Note that indeed
$X$ is a spine of $D^4$ that collapses to a point.
\end{example}

When the ribbon surface $S$ is a disc, it collapses to a point, and hence one might wonder whether we could simply take $X=S$ as a shadow. We show that this works only in the trivial case (the trivial ribbon disc is the one with one minimum and no saddle, having the trivial knot as boundary).

\begin{prop} \label{trivial:prop}
A properly embedded disc $D\subset D^4$ is a shadow if and only if $D$ is isotopic to the trivial ribbon disc (and hence $\partial D$ is the unknot).
\end{prop}
\begin{proof}
The disc $D$ collapses to a point, so $D$ is a shadow of $D^4$ if and only if $D^4 \searrow D$. This holds if and only if $D^4$ is a regular neighborhood of $D$. A regular neighborhood of $D$ is a product bundle $D\times D^2$, hence $D^4\searrow D$ if and only $D^4 =D\times D^2$. This holds precisely when $D$ is trivial.
\end{proof}

We have proved that every ribbon disc $D$ is contained in a shadow $X$, but $X$ may in fact be quite complicated.

\subsection{Non-ribbon surfaces} \label{non-ribbon:subsection}
One may wonder whether \emph{every} surface $S$ is contained in a shadow. We now show that this is not true: indeed being contained in a shadow is quite restrictive. Recall that a properly embedded surface $S\subset D^4$ is \emph{homotopically ribbon} if the inclusion 
$$(S^3\setminus \partial S) \hookrightarrow (D^4 \setminus S)$$
induces an epimorphism on fundamental groups
$$\pi_1(S^3\setminus \partial S) \twoheadrightarrow \pi_1(D^4 \setminus S).$$
For a general surface $S$, the following implications hold:
\begin{align} \label{implications:eqn}
S {\rm \ ribbon} \Longrightarrow S{\rm\ contained} {\rm\ in\ a\ shadow}\Longrightarrow S {\rm\ homotopically\ ribbon}.
\end{align}
We have already proved the first implication, so we now turn to the second. 
\begin{prop}
If $S$ is contained in a shadow $X$ then it is homotopically ribbon.
\end{prop}
\begin{proof}
The shadow $X$ contains $S$ and is hence obtained from $S$ by adding cells of index $0$, $1$, or $2$. Therefore a regular neighborhood $N(X)$ of $X$ is obtained from a regular neighborhood $N(S)$ of $S$ by adding handles of index $0$, $1$, or $2$. Since $X$ is a spine of $D^4$, we can take $N(X)=D^4$.

By turning handles upside-down we get that $D^4\setminus N(S)$ is obtained from a collar of $S^3\setminus N(\partial S)$ by adding handles of index $4$, $3$, or $2$. Since there are no 1-handles, the inclusion 
$$S^3\setminus N(\partial S) \hookrightarrow D^4 \setminus N(S)$$
induces a surjection on fundamental groups.
\end{proof}

We do not know if any of the two implications in (\ref{implications:eqn}) can be reversed. It is easy to construct some surface $S$ that is not homotopically ribbon, and such an $S$ cannot be contained in a shadow. The following example is certainly known to experts and we include it for completeness.

\begin{prop}
The trivial knot bounds some disc that is not homotopically ribbon.
\end{prop}
\begin{proof}
Pick a knotted sphere $S^2\subset S^4$ whose complement has non-cyclic fundamental group $G$, for instance a spun knot \cite[Chapter 3.J]{Rol}.

By tubing one such knotted sphere with a trivial properly embedded disc we get a disc $D^2\subset D^4$ such that $\pi_1(D^4\setminus D^2) = G$. Since $\partial D^2$ is the trivial knot, the complement $S^3\setminus \partial D^2$ is a solid torus and has cyclic $\pi_1$. The map 
$$\pi_1(S^3\setminus \partial D^2) \longrightarrow \pi_1(D^4 \setminus D^2)$$
cannot be surjective since the left group is cyclic and the right one is not.
\end{proof}

\subsection{Enlargement}

We prove here a stronger version of Theorem \ref{ribbon:contained:teo}:

\begin{teo} \label{ribbon:contained:2:teo}
Let $S\subset D^4$ be a ribbon surface and $G\subset S^3$ a knotted trivalent graph containing $\partial S$. There is a shadow $X$ of $G$ containing $S$.
\end{teo}
\begin{proof}
The ribbon surface $S$ is obtained from some planar diagram containing circles and edges as in Fig.~\ref{models:fig}-(left), and the link $L=\partial S\subset G$ is as in Fig.~\ref{models:fig}-(right). 

The graph $G$ contains $L$ and up to isotopy we may suppose that $G\setminus L$ is attached to $L$ only at the circles. Then we can proceed exactly as in the proof of Theorem \ref{ribbon:contained:teo} to get a shadow $X$ of $G$ containing $S$.
\end{proof}



\section{Shadows and the Kauffman bracket} \label{Jones:section}

We introduce the Kauffman bracket and Turaev's shadow formula.

\subsection{Kauffman bracket} \label{bracket:subsection}

\begin{figure}
\begin{center}
\includegraphics[width = 6 cm]{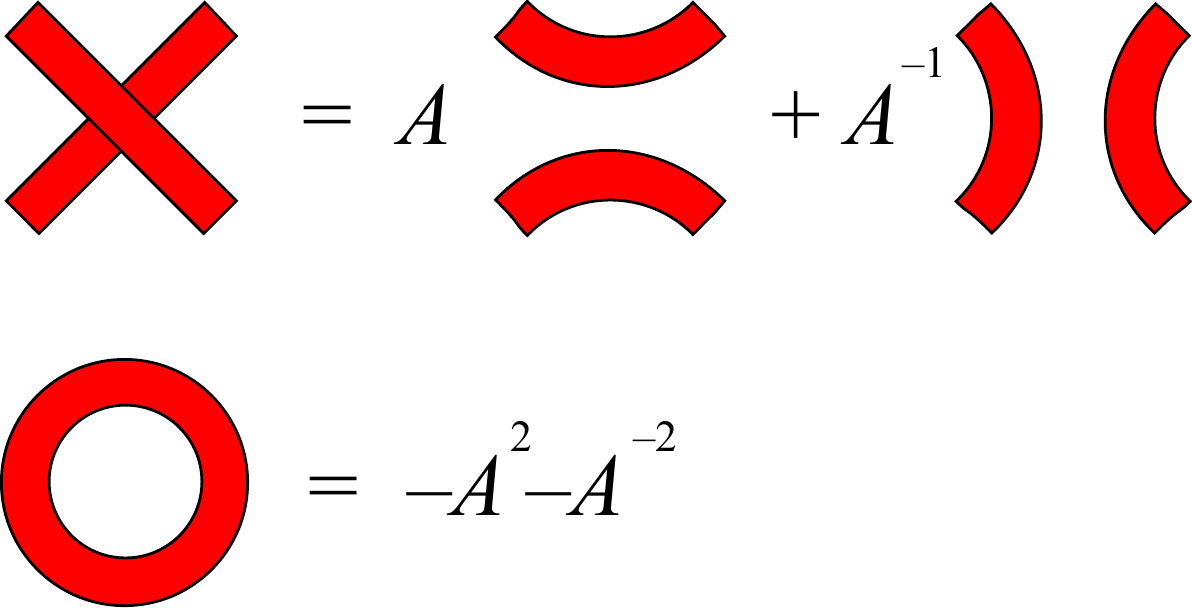}
\caption{The Kauffman bracket relations.}
\label{Kauffman_bracket:fig}
\end{center}
\end{figure}

The \emph{Kauffman bracket} $\langle L\rangle $ of a framed link $L\subset S^3$ is a polynomial in $\matZ[A,A^{-1}]$ defined using the skein relations shown in Fig.~\ref{Kauffman_bracket:fig}. The variables $q=A^2$ or $t=A^4$ are often used instead of $A$: the famous Jones polynomial of an oriented (but unframed) link is obtained from the Kauffman bracket simply by taking $t=A^4$ and assigning to the oriented link its Seifert framing. We will work with the variable $q = A^2$.

\subsection{Eisermann Theorem}

Eisermann has proved in \cite{Eis} the following fact.

\begin{teo} [Eisermann]
If $S\subset D^4$ is ribbon then $\langle\partial S\rangle$ has a zero in $q= i$ of order at least $\chi(S)$.
\end{teo}

The theorem provides some information only when $\chi(S)>0$. A $n$-component link $L$ is \emph{ribbon} if it bounds a ribbon surface consisting of $n$ discs. 

\begin{cor} \label{E:cor}
If a $n$-component link $L$ is ribbon then $\langle L \rangle$ has a zero at $q=i$ of order $n$.
\end{cor}

Eisermann has shown that this is the maximum order one can achieve: for every $n$-component link we have 
$$1\leqslant \ord_i \langle L \rangle \leqslant n$$
and both extremes may arise. In particular, when $n=1$ we always get $\ord_i \langle L \rangle = 1$ and hence Corollary \ref{E:cor} gives no information on knots.

Note that if we modify the framing of $L$ the Kauffman bracket $\langle L \rangle$ changes by a power of $A=q^{\frac 12}$ and hence its vanishing order at $q=i$ is unaffected: therefore we can neglect the framing in our investigation.

The Kauffman bracket of a link may also be calculated using shadows via a \emph{state-sum formula}. To explain this construction, due to Turaev, we need to introduce some objects.

\subsection{Colored ribbon graphs}
A \emph{framed knotted trivalent graph} $G\subset S^3$ is a knotted trivalent graph equipped with a \emph{framing}, \emph{i.e.}~an oriented surface thickening of the graph considered up to isotopy. An \emph{admissible coloring} of $G$ is the assignment of a natural number (a \emph{color}) at each edge of $G$ such that the three numbers $i,j,k$ coloring the three edges incident to a vertex satisfy the triangle inequalities, and their sum $i+j+k$ is even.

There is a standard way to define the Kauffman bracket $\langle G \rangle$ of a colored framed knotted trivalent graph $G\subset S^3$, which agrees with the above definition on framed links with all components colored by $1$. The bracket $\langle G \rangle$ will be a rational function in $q^{\frac 12}$ and not a Laurent polynomial in general -- although it turns out \emph{a posteriori} to be very close to a polynomial, see \cite{C:integral}.

\begin{figure}
 \begin{center}
   \includegraphics[width = 10 cm]{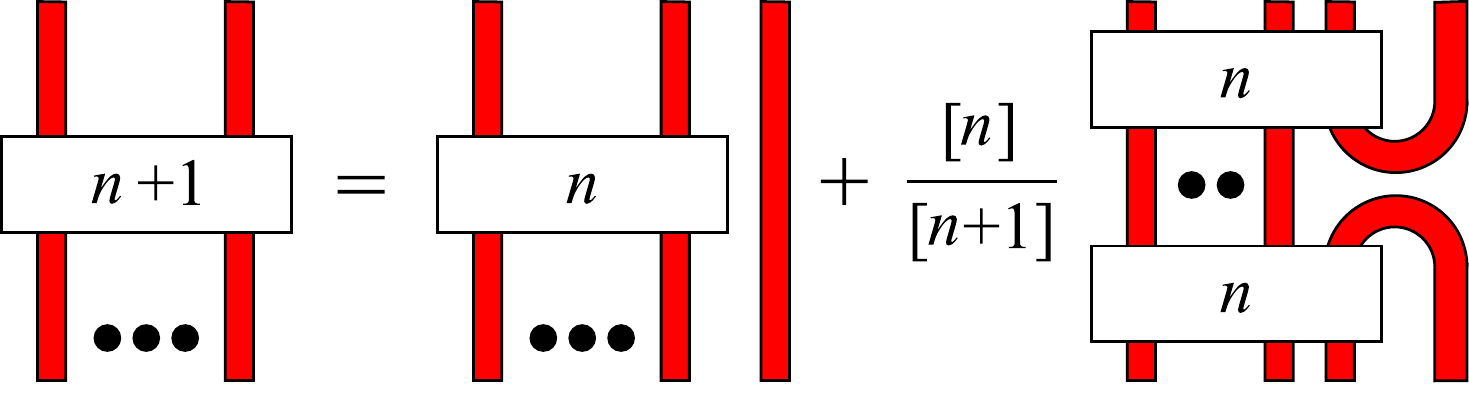}
   \end{center}
 \caption{The $(n+1)^{th}$ Jones-Wenzl projector is defined recursively with this formula.}
 \label{JW:fig}
\end{figure}

\begin{figure}
 \begin{center}
   \includegraphics[width = 8 cm]{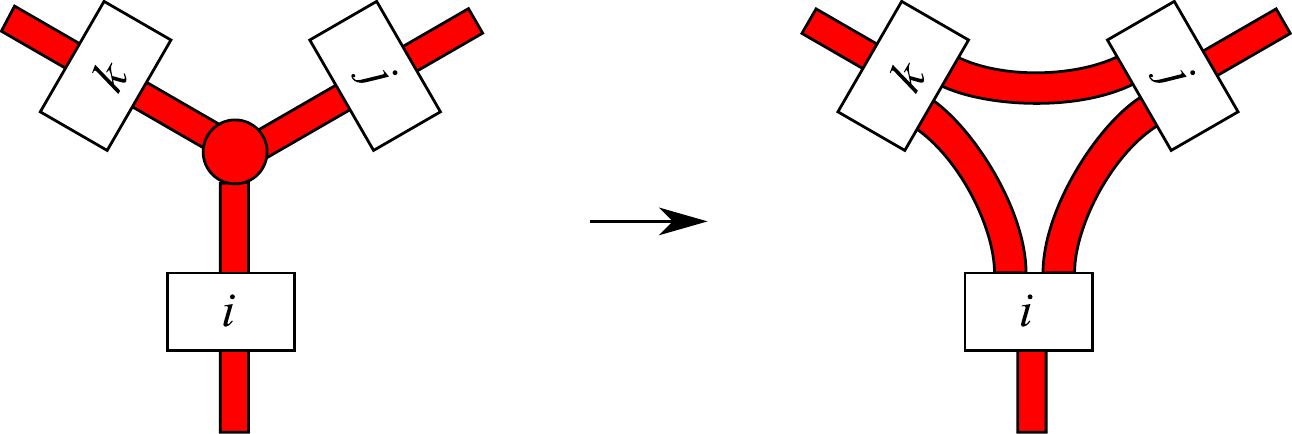}
   \end{center}
 \caption{A colored framed knotted trivalent graph determines a linear combination of framed links: replace every edge with a Jones-Wenzl projector, and connect them at every vertex via non intersecting strands contained in the depicted bands. For instance there are exactly $\frac{i+j-k}2$ bands connecting the projectors $i$ and $j$.}
 \label{ribbon_vertex:fig}
\end{figure}

To define $\langle G \rangle$ we must introduce the \emph{quantum integer}
$$[n] = \frac{q^{n}- q^{-n}}{q - q^{-1}} = q^{-n+1} + q^{-n+3} + \ldots + q^{n-3} + q^{n-1}.$$
The \emph{Jones-Wenzl projector} is a linear combination of framed arcs, defined recursively in Fig.~\ref{JW:fig}. The admissibility requirements on colors allow to associate uniquely to $G$ a linear combination of framed links by putting the $k^{\rm th}$ Jones-Wenzl projector at each edge colored with $k$ and by substituting vertices with bands as shown in Fig.~\ref{ribbon_vertex:fig}.

\subsection{Three important planar graphs} \label{basic:subsection}
\begin{figure}
 \begin{center}
  \includegraphics[width = 10.5 cm]{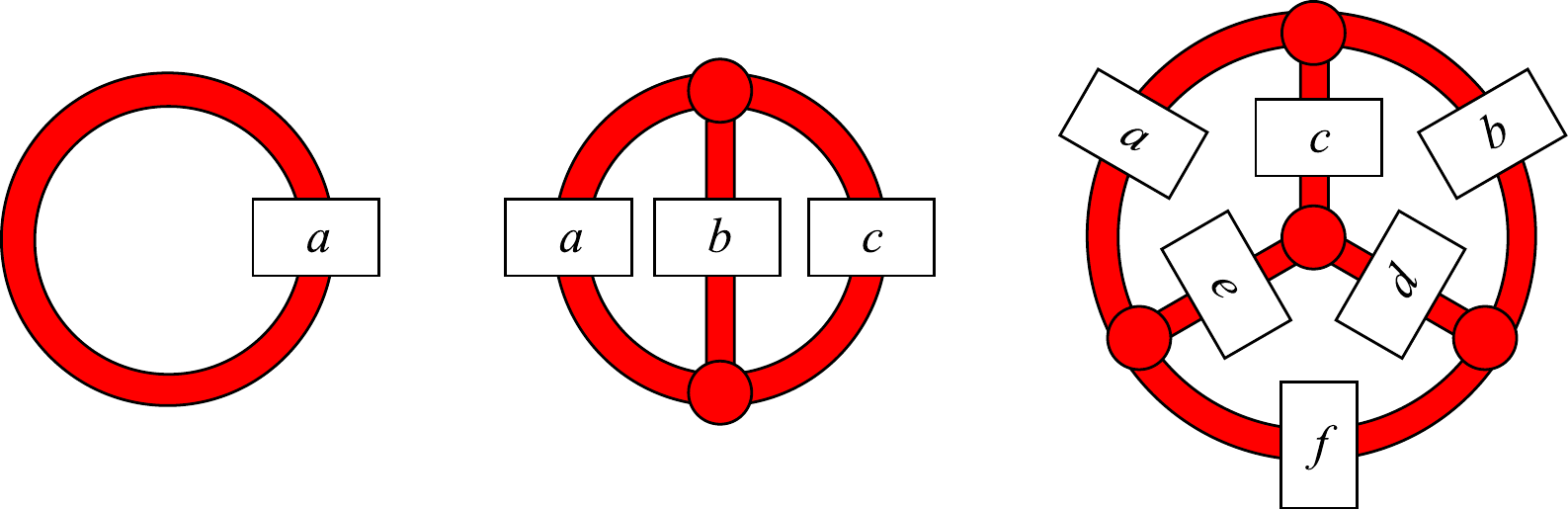}
 \end{center}
 \caption{Three important planar colored framed trivalent graphs in $S^3$.}
 \label{graphs:fig}
\end{figure}

Three basic planar framed trivalent graphs $\cerchio$, $\teta$, and $\tetra$ in $S^3$ are shown in Fig.~\ref{graphs:fig}. Their Kauffman brackets are some rational functions in $q$ that we now describe. 

We recall the usual factorial notation 
$$[n]! = [1]\cdots [n]$$ 
with the convention $[0]! = 1$. Similarly we define the generalized multinomials:
$$\begin{bmatrix} m_1, \ldots, m_h \\ n_1, \ldots n_k \end{bmatrix} = \frac{[m_1]!\cdots [m_h]!}{[n_1]!\cdots [n_k]!}.$$
When using these generalized multinomials we will always suppose that 
$$m_1 + \ldots + m_h = n_1+ \ldots + n_k.$$ 
The evaluations of $\cerchio$, $\teta$ and $\tetra$ are:
\begin{align*}
\cerchio_a  & = (-1)^{a} [a+1], \\
\teta_{a,b,c}  & = (-1)^{\frac{a+b+c}2} 
\begin{bmatrix} \frac{a+b+c}2+1, \frac{a+b-c}2, \frac{b+c-a}2, \frac{c+a-b}2 \\
a, b, c, 1 \end{bmatrix},  \\ 
\raisebox{-0.5cm}{\tetracolored } & = 
\begin{bmatrix} \Box_i-\triangle_j \\
a, b, c, d, e, f \end{bmatrix} \cdot \sum_{z = \max \triangle_j }^{\min \Box_i}\!\!\! (-1)^z 
\begin{bmatrix} z+1 \\
z-\triangle_j, \Box_i-z, 1 \end{bmatrix}.
\end{align*}

In the latter equality, triangles and squares are defined as follows:
\begin{align*}
\triangle_1 = \frac{a+b+c}{2},\ \triangle_2 = \frac{a+e+f}{2},\ \triangle_3 =\frac{ d+b+f}{2},\ \triangle_4 = \frac{d+e+c}{2},\\
\Box_1 = \frac{a+b+d+e}{2},\ \Box_2 = \frac{a+c+d+f}{2},\ \Box_3 = \frac{b+c+e+f}{2}.
\end{align*}

The indices in the formula vary as $1\leqslant i \leqslant 3$ and $1\leqslant j \leqslant 4$, so the term $\Box_i - \triangle_j$ indicates $3\times 4 = 12$ numbers.
The formula for $\tetra$ was first proved by Masbaum and Vogel \cite{MV}. These formulas are all rational functions in $q$ that may have poles in $0$, $\infty$, and at some root of unity, sometimes including the value $q=i$ we are interested in.

\subsection{Gleams}
Let $X\subset D^4$ be a shadow of a framed knotted trivalent graph $G\subset S^3$. Every region $R\subset X$ is equipped with a \emph{gleam}, a half-integer that generalizes the Euler number of closed surfaces embedded in oriented four-manifolds. The gleam is defined as follows.

\begin{figure}
 \begin{center}
  \includegraphics[width = 7 cm]{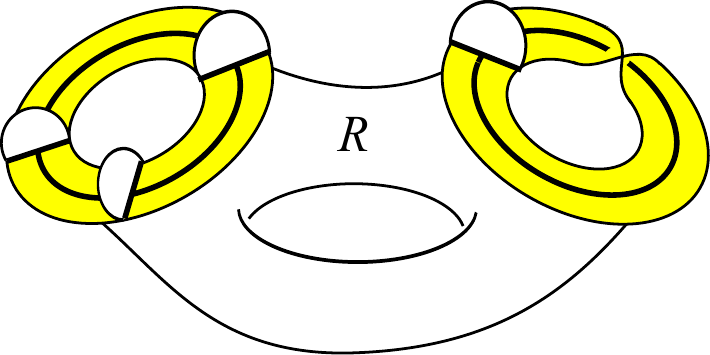}
 \end{center}
 \caption{A region $R$ in a shadow $X$. The shadow $X$ induces on every component of $\partial R$ an interval sub-bundle of the normal bundle in $D^4$, painted here in yellow.}
 \label{two-component:fig}
\end{figure}

The boundary $\partial R$ of $R$ consists of some closed curves, see Fig.~\ref{two-component:fig}. If $R$ is disjoint from $G$, the shadow $X$ provides an interval bundle over $\partial R$ as shown in the figure, which is an interval sub-bundle of the normal bundle of $\partial R$ in $D^4$. If $R$ is incident to some edge of $G$, the interval bundle is provided by the framing of $G$. (The boundary $\partial R$ is actually only immersed in general, but all these definitions work anyway.)

Let $R'$ be a generic small perturbation of $R$ with $\partial R'$ lying in the interval bundle at $\partial R$. The surfaces $R$ and $R'$ intersect only in isolated points, and we count them with signs:
$$\gl(R) = \frac 12 \# (\partial R \cap \partial R') + \#(R\cap R') \in \frac{\matZ}2$$
The half-integer $\gl(R)$ is the \emph{gleam} of $R$ and does not depend on the chosen $R'$. Note that the contribution of $\#(\partial R\cap \partial R')$ above one component of $\partial R$ is even or odd, depending on whether the interval bundle above it is an annulus or a M\"obius strip.

\subsection{Shadow formula}
Finally, we recall how to compute the Kauffman bracket of a colored framed knotted trivalent graph $G\subset S^3$ using shadows. 

Let $X$ be a shadow for $G$. An \emph{admissible coloring} $\sigma$ for $X$ is the assignment of a color to each region of $X$, such that for every interior edge of $X$ the colors of the three incident regions form an admissible triple. We also require that $\sigma$ extends the given coloring of $G$, \emph{i.e.}~a region $R\subset X$ incident to an edge $e$ of $G$ must be given the same color as $e$.

The \emph{evaluation} of the coloring $\sigma$ is the following function:
\begin{align} \label{evaluation:eqn}
\langle X_\sigma \rangle = \frac{\prod_f \cerchio_f^{\chi(f)}q_f \prod_v \tetra_v^{\chi(v)} \prod_{v_\partial} \teta_{v_\partial}^{\chi(v_\partial)}}
{\prod_e \teta_e^{\chi(e)} \prod_{e_\partial} \cerchio_{e_\partial}^{\chi(e_\partial)}}. 
\end{align}
Here the product is taken on all regions $f$, interior edges $e$, interior vertices $v$, boundary edges $e_\partial$, and boundary vertices $v_\partial$. The symbols
$$\cerchio_f, \teta_e, \tetra_v, \teta_{v_\partial}, \cerchio_{e_\partial}$$
indicate the Kauffman bracket of these graphs, colored respectively as $f$ or as the regions incident to $e, v, v_\partial, e_\partial$.

The \emph{phase} $q_f$ is the following monomial in $q^{\frac 14}$:
$$q_f = (\sqrt{-1})^{2gc}q^{-\frac{gc}2(c+2)}$$
where $g$ and $c$ are the gleam and the color of $f$, respectively.

The Euler characteristic $\chi(v)$ and $\chi(v_\partial)$ of vertices are obviously 1 and are included only for aesthetic reasons.

\begin{teo}[Turaev] \label{state-sum:teo}
Let $G\subset S^3$ be a colored framed knotted trivalent graph and $X$ any shadow for $G$. We have
$$\langle G \rangle = \sum_{\sigma} \langle X_\sigma \rangle$$
where the sum is taken on all colorings $\sigma$ of $X$ that extend that of $G$.
\end{teo}
We give a complete proof of this formula in Section \ref{state-sum:section}.

\section{Estimates at $q=i$} \label{estimate:section}
We prove all the needed estimates at $q=i$. The main result of this section is Theorem \ref{order:teo}, which is the technical core of the paper.

\subsection{Subsurfaces}
We will need the following.

\begin{prop} \label{correspondences:prop}
Let $X$ be a shadow of a trivalent knotted graph $G\subset S^3$. There are natural 1-1 correspondences:
$$
\left\{\!\begin{array}{c} {\rm properly\ embedded} \\ {\rm surfaces\ }S \subset X \end{array}\!\right\}
\longleftrightarrow
H_2(X,G;\matZ_2) 
\longleftrightarrow
H_1(G;\matZ_2) 
\longleftrightarrow 
\left\{\!\begin{array}{c} {\rm links} \\ L \subset G\end{array}\!\right\}.
$$
The correspondence sends $S$ to $L=\partial S$. The empty surface is included.
\end{prop}
\begin{proof}
The morphism $\partial\colon H_2(X,G;\matZ_2) \to H_1(G;\matZ_2)$ is an isomorphism because
$X$ is contractible and hence $H_i(X;\matZ_2)=\{e\}$ for $i=1,2$. Using cellular homology, every $\matZ_2$-homology class in $(X,G)$ is realized by a unique cycle, and that cycle is a surface since $X$ has simple singularities.
\end{proof}

Let now $\sigma$ be an admissible coloring for $X$. Its reduction modulo $2$ is a cycle in $H_2(X, G;\matZ_2)$ because the admissibility relation around every interior edge of $X$ reduces to $i+j+k \equiv 0$ (mod $2$). This cycle gives a surface $S_\sigma\subset X$ that consists of all regions in $X$ having an odd color: we call $S_\sigma$ the \emph{odd surface} of $\sigma$.

Analogously, an admissible coloring for $G$ determines an \emph{odd link} $L\subset G$ consisting of all edges with odd colors. Proposition \ref{correspondences:prop} implies the following:

\begin{cor}
Let $G\subset S^3$ be a colored framed knotted trivalent graph and $X$ be any shadow for $G$. The odd surface $S_\sigma\subset X$ of a coloring $\sigma$ that extends that of $G$ is the unique surface whose boundary $\partial S_\sigma$ is the odd sublink of $G$. In particular $S_\sigma$ does not depend on $\sigma$.
\end{cor}

\subsection{Red vertices}
Let $(a,b,c)$ be an admissible triple. Consider the following integers:
\begin{align} \label{triple:eqn}
\frac{a+b-c}2, \frac{b+c-a}2, \frac{c+a-b}2.
\end{align}

All the definitions we introduce are standard, except the following one which is new. We say that the triple $(a,b,c)$ is \emph{red} if at least two of the three integers in (\ref{triple:eqn}) are odd numbers. 

\begin{defn}
Let $G$ be a colored framed knotted trivalent graph. A vertex is \emph{red} if the colors of the three incident edges form a red triple.
\end{defn}

\subsection{The main technical theorem}

Given a meromorphic function $f(q)$ defined in a neighborhood of $q_0$, we denote by 
$$\ord_{q_0}f\in \matZ \cup\{+\infty\}$$ 
the maximum integer $k$ such that $f(q)/(q-q_0)^{k-1}$ vanishes in $q_0$. If $\ord_{q_0}f = + \infty$ the function $f$ vanishes in a neighborhood of $q_0$, otherwise it has a Laurent expansion
$$f(q) =\lambda (q-q_0)^{\ord_{q_0}f} + {\rm o}\big((q-q_0)^{\ord_{q_0}f}\big)$$
for some $\lambda \neq 0$.
We will be interested in the case $q_0=i$. We want to prove the following:

\begin{teo} \label{order:teo}
Let $X$ be a shadow colored by $\sigma$. We have 
$$\ord_i \langle X_\sigma \rangle \geqslant \chi(S_\sigma) - \frac{r}2$$ 
where $r$ is the number of red vertices in $\partial X$.
\end{teo}

This theorem and the topological Theorem \ref{ribbon:contained:teo} form altogether the core of this paper. In contrast with the topological one, this theorem has a long technical proof, to which we devote the rest of this section.

Before starting with the proof we single out some corollaries.

\begin{cor}
Let $G\subset S^3$ be a colored framed knotted trivalent graph. If the odd link $L\subset G$ bounds a ribbon surface $S\subset D^4$ then
$$\ord_i \langle G \rangle  \geqslant \chi(S) - \frac r2 $$ 
where $r$ is the number of red vertices in $G$.
\end{cor}
\begin{proof}
The ribbon surface $S$ is contained in a shadow $X$ of $G$ by Theorem \ref{ribbon:contained:2:teo}. We have
$$\langle G \rangle = \sum_\sigma \langle X_\sigma \rangle$$
which implies
$$\ord_i \langle G \rangle \geqslant \min_\sigma \big\{\ord_i \langle X_\sigma \rangle\big\}.$$
Every coloring $\sigma$ of $X$ extends the one of $G$ and hence its odd surface $S_\sigma \subset X$ has boundary $\partial S_\sigma = L$. Such a surface is unique by Proposition \ref{correspondences:prop} and hence necessarily $S_\sigma = S$. Now Theorem \ref{order:teo} says that
$$\ord_i \langle X_\sigma \rangle \geqslant \chi(S) - \frac r2 $$
for all $\sigma$. 
\end{proof}

A coloring of a link is \emph{odd} if each component is colored with an odd number.

\begin{cor}
If a link $L\subset S^3$ bounds a ribbon surface $S$ then
$$\ord_i \langle L \rangle \geqslant \chi(S)$$
for any framing and any odd coloring on $L$.
\end{cor}

Eisermann's Theorem corresponds to the case where all colorings are 1. 

\begin{cor}
Let $G$ be a colored framed knotted graph. If the odd link $L\subset G$ is ribbon, then
$$\ord_i\langle G \rangle \geqslant |L|- \frac r2 $$
where $|L|$ denotes the number of components of $L$ and $r$ is the number of red vertices in $G$.
\end{cor}
\begin{proof}
By hypothesis $L$ bounds a ribbon surface $S$ consisting with $|L|$ discs and hence $\chi(S) = |L|$.
\end{proof}

The following case includes the graphs \cerchio, \teta, and \tetra:

\begin{cor} \label{planar:cor}
Let $G\subset \matR^2$ be a colored planar graph. We have
$$\ord_i\langle G \rangle \geqslant |L|- \frac r2 $$
where $|L|$ denotes the number of components of the odd (un-)link $L\subset G$ and $r$ is the number of red vertices in $G$.
\end{cor}
\begin{proof}
The odd link $L$ is planar, hence trivial, hence ribbon.
\end{proof}

\subsection{Localization of Theorem \ref{order:teo}}
We now localize the proof of Theorem \ref{order:teo}, by reducing it to the building blocks \cerchio, \teta, and \tetra. The following lemma will be proved in the next section.

\begin{lemma} \label{block:lemma}
Let $G$ be a colored $\cerchio, \teta$, or $\tetra$. 
We have 
$$\ord_i\langle G \rangle \geqslant |L|- \frac r2 $$
where $L$ is the odd (un-)link $L\subset G$ and $r$  is the number of red vertices in $G$. If $G=\cerchio$ or $\teta$ then the equality holds.
\end{lemma}

Note that for $G=\cerchio, \teta, \tetra$ we have:
\begin{itemize}
\item $|L|=1$ if $G$ contains some odd-colored edges,
\item $|L|=0$ otherwise.
\end{itemize}
We postpone the proof of Lemma \ref{block:lemma} to the next section, and we now deduce Theorem \ref{order:teo} from it.

\vspace{2pt}\noindent\emph{Proof of Theorem \ref{order:teo}~from Lemma \ref{block:lemma}.} We have
$$\langle X_\sigma \rangle = \frac{\prod_f \cerchio_f^{\chi(f)}q_f \prod_v \tetra_v \prod_{v_\partial} \teta_{v_\partial}}
{\prod_e \teta_e^{\chi(e)} \prod_{e_\partial} \cerchio_{e_\partial}^{\chi(e_\partial)}}
$$
The phase $q_f$ is a monomial in $q$ and hence does not contribute to $\ord_i\langle X_\sigma \rangle$. We get
\begin{align*}
\ord_i\langle X_\sigma \rangle = & \sum_f \chi(f)\cdot \ord_i\cerchio_f + \sum_v \ord_i\tetra_v +  \sum_{v_\partial} \ord_i \teta_{v_\partial} \\ 
& - \sum_e \chi(e)\cdot \ord_i \teta_e - \sum_{e_\partial} \chi(e_\partial)\cdot\ord_i\cerchio_{e_\partial}.
\end{align*}
We now use Lemma \ref{block:lemma}. Note that for every colored $\cerchio, \teta, \tetra$ involved, we have $|L|=1$ precisely when the corresponding stratum (vertex, edge, or region) is contained in $S_\sigma$, otherwise we get $|L|=0$. We denote by $r(G)$ the number of red vertices in $G$ and we get:
\begin{align*}
\ord_i\langle X_\sigma \rangle \geqslant & \sum_{f\subset S_\sigma} \chi(f) + \sum_{v\in S_\sigma} 1 +  \sum_{v_\partial\in S_\sigma} 1 
- \sum_{e\subset S_\sigma} \chi(e) - \sum_{e_\partial\subset S_\sigma} \chi(e_\partial) \\
& - \sum_v \frac {r(v)}2  
- \sum_{v_\partial}\frac {r(v_\partial)}2 
+ \sum_e \frac {r(e)}2 \\
= & \ \chi(S_\sigma) - \sum_v \frac {r(v)}2 
- \sum_{v_\partial}\frac {r(v_\partial)}2
+ \sum_e \frac {r(e)}2.
\end{align*}
Let $e$ be an interior edge. The two vertices of $\teta_e$ are colored by the same triple $(a,b,c)$: hence $\teta_e$ has either zero or two red vertices. If an interior vertex $v$ of $X$ is adjacent to $e$, then $\tetra_v$ has a corresponding vertex colored by $(a,b,c)$. If an exterior vertex $v_\partial$ is adjacent to $e$, then both vertices of $\teta_{v_\partial}$ are colored as $(a,b,c)$. From this we get
$$\sum_e r(e) = \sum_v r(v) + \sum_{v_\partial} \frac{r(v_\partial)}2$$
and therefore
$$ \ord_i\langle X_\sigma \rangle \geqslant \chi(S_\sigma) - \sum_{v_\partial} \frac {r(v_\partial)}4 = \chi(S_\sigma) - \frac r2$$
because $r(v_\partial)$ equals $2$ when $v_\partial$ is red and $0$ otherwise.
\finedimo

\subsection{Order of generalized multinomials}
It remains to prove Lemma \ref{block:lemma}, and to do so we will need the following.
\begin{prop} \label{orders:prop}
We have
\begin{align*}
\ord_i[n] & = 
\left\{\begin{array}{ll} 0 & {\rm if\ } n\ {\rm is\ odd},  \\
1 & {\rm if\ } n\ {\rm is\ even}, 
\end{array} \right. \\
\ord_i[n]! & = \bigg\lfloor \frac n2 \bigg\rfloor, \\
\ord_i \begin{bmatrix} m_1, \ldots, m_h \\ n_1, \ldots n_k \end{bmatrix} & = 
\Bigg\lfloor \frac{\#\big\{{\rm odd}\ n_i \big\}}2\Bigg\rfloor - \Bigg\lfloor\frac{\#\big\{{\rm odd\ } m_j \big\}}2 \Bigg\rfloor.
\end{align*}
\end{prop}
\begin{proof}
The function
$$[n] = \frac {q^n - q^{-n}}{q-q^{-1}} = \frac{q^{-n}}{q-q^{-1}} (q^{2n}-1)$$
has simple zeroes at the $(2n)^{\rm th}$ roots of unity (except $q=\pm 1$), hence at $q=i$ when $n$ is even. The equality $\ord_i[n]! = \lfloor \tfrac n2 \rfloor$ follows. On the multinomial, recall that $m_1+\ldots +m_h = n_1+\ldots+n_h = N$ by hypothesis. We get
\begin{align*}
\ord_i \begin{bmatrix} m_1, \ldots, m_h \\ n_1, \ldots n_k \end{bmatrix} 
& = \sum_i \bigg\lfloor \frac{m_i}2 \bigg\rfloor - \sum_j \bigg\lfloor \frac{n_j}2 \bigg\rfloor
\\
& = \bigg\lfloor\frac N2\bigg\rfloor - \bigg\lfloor\frac{\#\big\{{\rm odd\ } m_i\big\}}2\bigg\rfloor - 
\bigg\lfloor\frac N2\bigg\rfloor + \bigg\lfloor\frac{\#\big\{{\rm odd\ } n_j\big\}}2\bigg\rfloor \\
& = \Bigg\lfloor \frac{\#\big\{{\rm odd}\ n_i \big\}}2\Bigg\rfloor - \Bigg\lfloor\frac{\#\big\{{\rm odd\ } m_j \big\}}2 \Bigg\rfloor.
\end{align*}
\end{proof}

We can now evaluate $\cerchio, \teta$, and $\tetra$ at $q=i$. 

\subsection{Orders of the circle, theta, and tetrahedron}
It remains to prove Lemma \ref{block:lemma}

\dimo{block:lemma}
If $G=\cerchio$ then
$\ord_i\cerchio_a = \ord_i[a+1]$ equals $1$ if $a$ is odd and $0$ if $a$ is even: the odd link $L$ is respectively $G$ and $\emptyset$, therefore $\ord_i\cerchio_a = |L|$ in any case.

If $G=\teta$ we have
\begin{align*}
\ord_i\teta_{a,b,c} & = \ord_i \begin{bmatrix} \frac{a+b+c}2+1, \frac{a+b-c}2, \frac{b+c-a}2, \frac{c+a-b}2 \\
a, b, c, 1 \end{bmatrix} \\
 & = \Bigg\lfloor \frac{\#\big\{{\rm odd}\ a,b,c,1 \big\}}2 \Bigg\rfloor - 
 \Bigg\lfloor\frac{\#\big\{{\rm odd\ } \frac{a+b+c}2+1, \frac{a+b-c}2, \frac{b+c-a}2, \frac{c+a-b}2 \big\}}2 \Bigg\rfloor \\
 & = |L| - \frac r2.
\end{align*}
To prove the last equality, note that the first addendum is $0$ if $a,b,c$ are even and $1$ otherwise (there are either zero or two odd numbers in $a,b,c$ by admissibility), and $L\subset G$ is respectively empty or a circle.
Concerning the second addendum, note that
$$\frac{a+b+c}2 +1= \frac{a+b-c}2 + \frac{b+c-a}2 + \frac{c+a-b}2 +1$$
and hence one easily sees that the second addendum equals
$$\Bigg\lfloor\frac{\#\big\{{\rm odd\ } \frac{a+b-c}2, \frac{b+c-a}2, \frac{c+a-b}2 \big\}}2 \Bigg\rfloor$$
which is 1 if the triple is red and 0 otherwise, by definition.

For $G=\tetra$ we do a long case-by-case analysis. We recall the formula
\begin{align*}
\raisebox{-0.5cm}{\tetracolored } & = 
\begin{bmatrix} \Box_i-\triangle_j \\
a, b, c, d, e, f \end{bmatrix} \cdot \sum_{z = \max \triangle_j }^{\min \Box_i}\!\!\! (-1)^z 
\begin{bmatrix} z+1 \\
z-\triangle_j, \Box_i-z, 1 \end{bmatrix}.
\end{align*}
with 
\begin{align*}
\triangle_1 = \frac{a+b+c}{2},\ \triangle_2 = \frac{a+e+f}{2},\ \triangle_3 =\frac{ d+b+f}{2},\ \triangle_4 = \frac{d+e+c}{2},\\
\Box_1 = \frac{a+b+d+e}{2},\ \Box_2 = \frac{a+c+d+f}{2},\ \Box_3 = \frac{b+c+e+f}{2}.
\end{align*}
Note that
$$a+b+c+d+e+f = \sum_i \Box_i = \sum_j \triangle_j.$$
We now estimate the factor
\begin{equation} \label{factor2:eqn}
\sum_{z = \max \triangle_j }^{\min \Box_i}\!\!\! (-1)^z 
\begin{bmatrix} z+1 \\
z-\triangle_j, \Box_i-z, 1 \end{bmatrix}
\end{equation}
in terms of the parity of the $\Box_j$'s and the $\triangle_i$'s.

We first consider the case $a+b+c+d+e+f$ is even. In that case the number of odd $\Box_i$'s is 0 or 2, while the number of odd $\triangle_j$'s is 0, 2, or 4. Using Proposition \ref{orders:prop} we easily see that
$$\ord_i \begin{bmatrix} z+1 \\
z-\triangle_j, \Box_i-z, 1 \end{bmatrix}$$
is a number that depends on the parity of $z$, on the number $0,2$ of odd $\Box_i$'s and $0,2,4$ of odd $\triangle_j$'s according to the tables: 
\begin{center}
\begin{tabular}{c|c|c}
\multicolumn{3}{c}{$z$ even} \\
\hline \hline
 & $0\ \Box_i$ & $2\ \Box_i$  \\
\hline
 $0\ \triangle_j$ & $0$ & $1$ \\
 $2\ \triangle_j$ & $1$ & $2$ \\
 $4\ \triangle_j$ & $2$ & $3$
 \end{tabular}
\qquad
\begin{tabular}{c|c|c}
\multicolumn{3}{c}{$z$ odd} \\
\hline \hline
 & $0\ \Box_i$ & $2\ \Box_i$  \\
\hline
 $0\ \triangle_j$ & $4$ & $3$ \\
 $2\ \triangle_j$ & $3$ & $2$ \\
 $4\ \triangle_j$ & $2$ & $1$
 \end{tabular}
\end{center}
By taking the minimum we get that the order at $q=i$ of (\ref{factor2:eqn}) is at least:
\begin{equation} \label{even:eqn}
\begin{tabular}{c|c|c} 
 & $0\ \Box_i$ & $2\ \Box_i$  \\
\hline
 $0\ \triangle_j$ & $0$ & $1$ \\
 $2\ \triangle_j$ & $1$ & $2$ \\
 $4\ \triangle_j$ & $2$ & $1$
 \end{tabular}
\end{equation}
The case $a+b+c+d+e+f$ odd is treated analogously: now the number of odd $\Box_i$'s is $1$ or $3$, and the number of odd $\triangle_i$'s is $1$ or $3$. We get
\begin{center}
\begin{tabular}{c|c|c}
\multicolumn{3}{c}{$z$ even} \\
\hline \hline
 & $1\ \Box_i$ & $3\ \Box_i$  \\
\hline
 $1\ \triangle_j$ & $1$ & $2$ \\
 $3\ \triangle_j$ & $2$ & $3$
 \end{tabular}
\qquad
\begin{tabular}{c|c|c}
\multicolumn{3}{c}{$z$ odd} \\
\hline \hline
 & $1\ \Box_i$ & $3\ \Box_i$  \\
\hline
 $1\ \triangle_j$ & $3$ & $2$ \\
 $3\ \triangle_j$ & $2$ & $1$
 \end{tabular}
\end{center}
The order at $q=i$ of (\ref{factor2:eqn}) is hence at least:
\begin{equation} \label{odd:eqn}
\begin{tabular}{c|c|c} 
 & $1\ \Box_i$ & $3\ \Box_i$  \\
\hline
 $1\ \triangle_j$ & $1$ & $2$ \\
 $3\ \triangle_j$ & $2$ & $1$
 \end{tabular}
\end{equation}
We now turn to the factor
\begin{equation} \label{factor1:eqn}
\begin{bmatrix} \Box_i-\triangle_j \\
a, b, c, d, e, f \end{bmatrix}.
\end{equation}
The 12 numbers $\Box_i - \triangle_j$ are of type $\frac{m+n-p}2$ where $(m,n,p)$ are the colors of the edges incident to some vertex: there are $4$ vertices and $3$ such expressions at each vertex; the $12$ numbers correspond to the 12 red arcs in the picture
\begin{center}
\includegraphics[width = 1.5 cm]{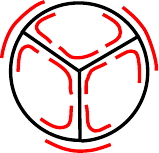}
\end{center}
where the red arc corresponding to $\frac{m+n-p}2$ is the one parallel to the edges $m,n$ and opposite to $p$. The parities of these 12 numbers determine the parities of all the quantities $\Box_i, \triangle_j, a,b,c,d,e,f$, and hence also $|L|$ and $\frac r2$.
The possible configurations (considered up to symmetries of the tetrahedron) are easily classified and are shown in Tables \ref{cases:table} and \ref{cases2:table}.

\begin{table}
\begin{center}
\begin{tabular}{c|c|c|c|c|c|c|c}
\phantom{\Big|}\! odd $\Box_i$'s & odd $\triangle_j$'s & red arcs & $\ord_i\big((\ref{factor1:eqn})\big)$ & $\ord_i\big((\ref{factor2:eqn})\big)$ & $|L|$ & $\frac r2$ & works? \\
\hline \hline
 0 & 0 & \raisebox{-0.65 cm}{\includegraphics[width = 1.8 cm]{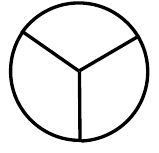}} & 0 & $\geqslant 0$ & $0$ & $0$ & yes \\
 0 & 2 & \raisebox{-0.65 cm}{\includegraphics[width = 1.8 cm]{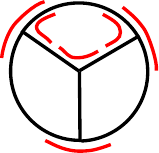}} & $-1$ & $\geqslant 1$ & $1$  & 1& yes \\
  0 & 4 & \raisebox{-0.65 cm}{\includegraphics[width = 1.8 cm]{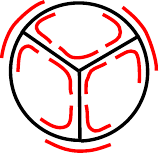}} & $-6$ & $\geqslant 2$ & $0$ & 2& no \\
  2 & 0 & \raisebox{-0.65 cm}{\includegraphics[width = 1.8 cm]{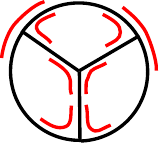}} & $-2$ & $\geqslant 1$ & $1$  & 2& yes \\
  2 & 2 & \raisebox{-0.65 cm}{\includegraphics[width = 1.8 cm]{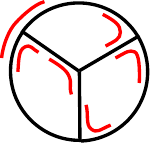}} & $-1$ & $\geqslant 2$ & $1$  & 1& yes \\
  2 & 2 & \raisebox{-0.65 cm}{\includegraphics[width = 1.8 cm]{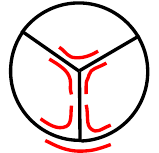}} & $-3$ & $\geqslant 2$ & $0$  & 1& yes \\
  2 & 4 & \raisebox{-0.65 cm}{\includegraphics[width = 1.8 cm]{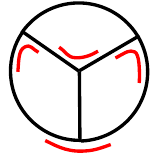}} & $0$ & $\geqslant 1$ & $1$  & 0 & yes 
\end{tabular}
\caption{For each case: the number of odd $\Box_i$'s, of odd $\triangle_j$'s, the red arcs, the order of the first factor (\ref{factor1:eqn}), of the second (\ref{factor2:eqn}) estimated in (\ref{even:eqn}), the number of components of the odd link $|L|$, and $\frac r2$. If 
(\ref{factor1:eqn}) + (\ref{factor2:eqn}) $\geqslant |L|-\frac r2$ then the estimate works (last column).}
\label{cases:table}
\end{center}
\end{table}

\begin{table}
\begin{center}
\begin{tabular}{c|c|c|c|c|c|c|c}
\phantom{\Big|}\! odd $\Box_i$'s & odd $\triangle_j$'s & red arcs & $\ord_i\big((\ref{factor1:eqn})\big)$ & $\ord_i\big((\ref{factor2:eqn})\big)$ & $|L|$ & $\frac r2$ & works? \\
\hline \hline
 1 & 1 & \raisebox{-0.65 cm}{\includegraphics[width = 1.8 cm]{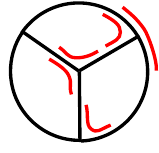}} & $-1$ & $\geqslant 1$ & $1$ & $1$ & yes \\
 1 & 3 & \raisebox{-0.65 cm}{\includegraphics[width = 1.8 cm]{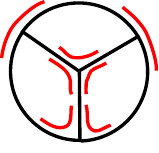}} & $-2$ & $\geqslant 2$ & $1$ & $1$ & yes \\
 3 & 1 & \raisebox{-0.65 cm}{\includegraphics[width = 1.8 cm]{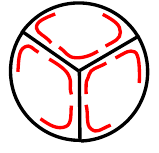}} & $-3$ & $\geqslant 2$ & $1$ & $2$ & yes \\
 3 & 3 & \raisebox{-0.65 cm}{\includegraphics[width = 1.8 cm]{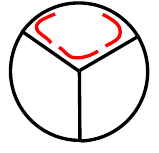}} & 0 & $\geqslant 1$ & $1$ & $0$ & yes
\end{tabular}
\caption{For each case: the number of odd $\Box_i$'s, of odd $\triangle_j$'s, the red arcs, the order of the first factor (\ref{factor1:eqn}), of the second (\ref{factor2:eqn}) estimated in (\ref{odd:eqn}), the number of components of the odd link $|L|$, and $\frac r2$. If 
(\ref{factor1:eqn}) + (\ref{factor2:eqn}) $\geqslant |L|-\frac r2$ then the estimate works (last column).}
\label{cases2:table}
\end{center}
\end{table}

As the tables show, the needed inequality 
$$\ord_i\big((\ref{factor2:eqn})\big) + \ord_i\big((\ref{factor1:eqn})\big) \geqslant |L| + \frac r2$$
is verified for all the configurations, except one bad case: when the $\Box_i$'s are all even and the $\triangle_j$'s are all odd we need to prove that 
$$\ord_i\big((\ref{factor2:eqn})\big) + \ord_i\big((\ref{factor1:eqn})\big) \geqslant -2$$
but we only get $\geqslant -4$. This bad case holds for instance when $a=b=c=d=e=f=2$ and hence $\Box_i = 4$ and $\triangle_j = 3$. If we look more carefully at this example we find

\begin{align*}
\raisebox{-0.5cm}{\tetracoloredtwo} & = 
\begin{bmatrix} 1 \cdots 1 \\
2, 2, 2, 2, 2, 2 \end{bmatrix} \cdot \sum_{z = 3}^{4} (-1)^z 
\begin{bmatrix} z+1 \\
z-3, 4-z, 1 \end{bmatrix} \\
& = \frac{1}{[2]^6} \cdot \left(-[4]!+[5]!\right) \\
& = \frac{[4]!}{[2]^6} \cdot ([5]-1).
\end{align*}
Now it turns out that the difference
$$[5]-1 = q^4+q^2+q^{-2}+q^4 = (q+q^{-1})(q^3+q^{-3}) = [2]\big([4]-[2]\big)$$
has order $2$ at $q=i$: this difference produces a cancellation that increases the order of (\ref{factor2:eqn}) at $q=i$ by two, giving overall the desired $-4$ instead of the $\geqslant -2$ expected by the tables.

We now prove that this kind of cancelation holds in general, provided that the $\Box_i$'s are all even and the $\triangle_j$'s are all odd. The sum
\begin{equation*}
\sum_{z = \max \triangle_j }^{\min \Box_i}\!\!\! (-1)^z 
\begin{bmatrix} z+1 \\
z-\triangle_j, \Box_i-z, 1 \end{bmatrix}
\end{equation*}
goes from the odd $z=\max \triangle_j$ to the even $z=\min \Box_i$ and so contains an even number of terms. Two subsequent terms $z=2k-1$ and $z=2k$ give
$$
- \begin{bmatrix} 2k \\
2k-1-\triangle_j, \Box_i-2k+1, 1 \end{bmatrix}
+ \begin{bmatrix} 2k+1 \\
2k-\triangle_j, \Box_i-2k, 1 \end{bmatrix}
$$
that may be rewritten as
$$
- \begin{bmatrix} 2k \\ 2k-1-\triangle_j, \Box_i-2k,1,1,1,1 \end{bmatrix} \cdot \left( \frac 1{\prod_i[\Box_i-2k+1]} - \frac{[2k+1]}{\prod_j[2k-\triangle_j]}\right).$$
The left factor has order $2$ as prescribed by Table \ref{cases:table}. Quite surprisingly, the second factor 
$$\frac{\prod_j[2k-\triangle_j] - {[2k+1]}\cdot \prod_i[\Box_i-2k+1]}
{\prod_i[\Box_i-2k+1]\cdot \prod_j[2k-\triangle_j]}.$$
has order at least $2$: note that all the quantum integers in the formula are quantum odd numbers; the denominator is a non-zero constant at $q=i$, while the numerator has order $\geqslant 2$ thanks to the following lemma. 

\begin{lemma} \label{bad:case:lemma}
Let $x_1,\ldots,x_n,y_1,\ldots,y_m$ be odd non-negative integers with
$$\sum_j (y_j-1) \equiv \sum_i (x_i-1) \ ({\rm mod}\ 4).$$
Then
$$\ord_i\left(\prod_i[x_i] - \prod_j[y_j]\right) \geqslant 2.$$
\end{lemma}
\begin{proof}
We set $f(q) = \prod_i[x_i] - \prod_j[y_j]$ and write $\sqrt{-1}$ instead of $i$ to avoid confusion. Now
$$[2k+1](\sqrt{-1}) = (-1)^k$$
gives
$$ f(\sqrt{-1})  = (-1)^{\frac 12\sum_i (x_i-1)} - (-1)^{\frac 12 \sum_j (y_j -1)} = 0$$
since $\frac 12\sum_i (x_i-1)$ and $\frac 12 \sum_j (y_j -1)$ have the same parity by hypothesis. This gives $\ord_if \geqslant 1$. We now calculate the derivative $f'$ of $f$. Note that
$$[n]' = \frac{n(q^{n-1} + q^{-n-1})(q-q^{-1}) - (1+q^{-2})(q^n-q^{-n})}{(q-q^{-1})^2}.$$
vanishes when $q=\sqrt{-1}$ and $n$ is odd, since both $q^{n-1}+q^{-n-1}$ and $1+q^{-2}$ do. Therefore the derivatives of $\prod [x_i]$ and $\prod [y_j]$ both vanish at $q=\sqrt{-1}$ and hence $f'(\sqrt{-1})=0$. Therefore $\ord_i f \geqslant 2$.
\end{proof}

To conclude the proof of Lemma \ref{block:lemma} we must verify that 
\begin{equation*} 
\sum_j (2k-\triangle_j-1) \equiv 2k + \sum_i (\Box_i-2k) \ ({\rm mod}\ 4)
\end{equation*}
and apply Lemma \ref{bad:case:lemma}. This is equivalent to $\sum_j \triangle_j \equiv \sum_i \Box_i$ which is true since actually $\sum_j \triangle_j = \sum_i \Box_i$.
\finedimo

\section{Other manifolds} \label{other:section}
We extend everything from $S^3$ to $\#_g(S^2\times S^1)$. We prove in particular Theorems \ref{link:teo} and \ref{large:teo}.

\subsection{Ribbon surfaces}
The notion of ribbon surface extends naturally from $S^3$ to every closed 3-manifold $M$.
A properly embedded smooth surface $S\subset M\times [0,1]$ with $\partial S \subset M \times 0$ is \emph{ribbon} if one of the following equivalent conditions holds:

\begin{itemize}
\item the surface $S$ may be isotoped to an immersed surface in $M$ having only ``ribbon'' singularities as in Fig.~\ref{ribbon_singularity:fig},
\item the surface $S$ may be isotoped in $M\times [0,1]$ into Morse position, with only minima and saddle points (no maxima) as in Fig.~\ref{Morse:fig}.
\end{itemize}

Every ribbon surface $S$ can be constructed from a graph embedded in $M$ as in Fig.~\ref{construct_ribbon:fig}, consisting of some disjoint circles bounding discs (the minima), and some arcs connecting them in space (the saddles).

\subsection{Shadows}
Our definition of shadow is very restrictive and designed for $D^4$, and it cannot be extended harmlessly to manifolds other than $S^3$.

Costantino has proposed in \cite{C:Jones} a definition when $M_g=\#_g(S^2\times S^1)$ is a connected sum of some $g\geqslant 1$ copies of $S^2\times S^1$. In that case $M_g$ is the boundary of the oriented four-dimensional handlebody $H^4_g$ made of one 0-handle and $g$ one-handles.

\begin{defn}
A \emph{shadow} for $H^4_g$ is a simple polyhedron $X\subset H^4_g$ such that the following holds:
\begin{itemize}
\item $X$ is properly embedded, that is $\partial X = X\cap M_g$,
\item $X$ is locally flat,
\item $X$ collapses onto a graph $Y$,
\item $H^4_g$ collapses onto $X$.
\end{itemize}
\end{defn}
The last two conditions can be summarized by writing
$$H^4_g \searrow X \searrow Y.$$
The boundary $G=\partial X$ of a shadow $X\subset H_g^4$ is a knotted trivalent graph $G\subset M_g$, and we say that $X$ is a shadow of $G$.

\begin{prop}[Costantino \cite{C:Jones}] \label{Cost:prop}
Every knotted trivalent graph $G\subset M_g$ has a shadow $X\subset H_g^4$.
\end{prop}
\begin{proof}
Same proof as in Proposition \ref{Turaev:prop}, with a small variation. We set $H_g^4 = D_g\times D^2$ where $D_g$ is a disc with $g$ holes. Up to isotopy we can see $G$ as a diagram in the interior of $D_g$. Up to some Reidemeister move we can suppose that there is a smallest closed disc with $g$ holes $D_g'\subset D_g$ containing $G$ such that $D_g \setminus D_g'$ is a collar of $\partial D_g$. We clearly have $H_g^4 \searrow D_g' \searrow Y$ for some graph $Y\subset D_g'$.
We enlarge $D_g'$ by adding a cylinder above $G$ and we get a shadow $X$ for $G$.
\end{proof}

\subsection{Ribbon surfaces in a shadow}
We can now extend Theorem \ref{main:topological:teo} from $S^3$ to $M_g$.
A ribbon surface in a 4-manifold like $H_g^4$ is just a ribbon surface in a collar of its boundary.

\begin{teo} \label{large:topological:teo}
Every ribbon surface $S\subset H_g^4$ is contained in a shadow $X$ with $\partial X = S$.
\end{teo}
\begin{proof}
Same proof as in Theorem \ref{ribbon:contained:teo}. We construct $S$ from a graph $G\subset \#_g(S^2\times S^1)$ as in Fig.~\ref{construct_ribbon:fig} made of circles and arcs. Up to Reidemeister moves we suppose that $G$ is contained in a smallest disc with holes and we build a shadow $X$ for $G$ as described in the proof of Proposition \ref{Cost:prop}. Then we add bands and push them inside $H_g^4$.
\end{proof}

\subsection{Kauffman bracket}
The Kauffman bracket is also defined in $M_g$, thanks to result of Hoste-Przytycki \cite{HP, P2} and (with different techniques) to Costantino \cite{C:Jones}. We briefly recall its definition.

Let $M$ be an oriented 3-manifold. Consider the field $\matQ(A)$ of all complex rational functions with variable $A$ and the abstract $\matQ(A)$-vector space $V$ generated by all framed links in $M$, considered up to isotopy. The \emph{skein vector space} $K(M)$ is the quotient of $V$ by all the possible skein relations as in Fig.~\ref{Kauffman_bracket:fig}. An element of $K(M)$ is called a \emph{skein}.

\begin{prop}
The skein vector space $K(M_g)$ of $M_g$ is isomorphic to $\matQ(A)$ and generated by the empty skein $\emptyset$.
\end{prop}
\begin{proof}
This is due to Hoste and Przytycki \cite{HP, P, P2}, see also \cite[Proposition 1]{FK}.
\end{proof}

A colored framed knotted trivalent graph $G$ determines a skein $G\in K(M)$ and as such it is equivalent to $\langle G \rangle \cdot \emptyset$ for a unique coefficient $\langle G \rangle \in \matQ(A)$. This coefficient is by definition the Kauffman bracket $\langle G \rangle$ of $G$.

\begin{rem} \label{canonical:rem}
There is an obvious canonical linear map $K(M) \to K(M\# N)$ defined by considering a skein in $M$ inside $M\# N$. The linear map $K(M_g) \to K(M_{g+1})$ sends $\emptyset$ to $\emptyset$ and hence preserves the bracket $\langle G \rangle$ of a $G\subset M_g$.

This shows in particular that if $G$ is contained in a ball, the bracket $\langle G \rangle$ is the same that we would obtain by considering $G$ inside $S^3$.
\end{rem}

\subsection{Shadow formula}

The shadow formula works also in this context.

\begin{teo}[Shadow formula] \label{state-sum2:teo}
Let $G\subset M_g$ be a colored framed knotted trivalent graph and $X\subset H^4_g$ any shadow for $G$. We have
$$\langle G \rangle = \sum_{\sigma} \langle X_\sigma \rangle$$
where the sum is taken on all colorings $\sigma$ of $X$ that extend that of $G$.
\end{teo}

A crucial observation \cite[Lemma 3.6]{C:Jones} is that the number of colorings $\sigma$ extending that of $G$ is finite, because $X$ collapses to a graph $Y$: hence the sum makes sense (see Proposition \ref{collapses:prop}). We prove Theorem \ref{state-sum2:teo} in Section \ref{state-sum:section}.

\begin{rem}
Costantino \cite{C:Jones} uses the shadow formula to \emph{define} $\langle G \rangle$ and then employs Turaev's theory of shadows and Reshetikhin-Turaev invariants to prove that the result does not depend on the shadow chosen. Costantino's definition agrees with ours up to a slightly different normalization: he wants to extend the Jones polynomial, while we prefer to extend the Kauffman bracket. 
\end{rem}

\subsection{Main theorem}
We can finally prove Theorem \ref{large:teo}:

\begin{teo}  \label{large:2:teo}
Let $G$ be a colored framed knotted trivalent graph in $S^3$ or $\#_g(S^2\times S^1)$ and $L\subset G$ be its odd sublink. If $L$ bounds a ribbon surface $S$ then
$$\ord_i\langle G \rangle \geqslant \chi(S) - \frac r2$$
where $r$ is the number of red vertices in $G$.
\end{teo}
\begin{proof}
We know that $L$ has a shadow $X$ containing $S$ by Theorem \ref{large:topological:teo}. The proof of Theorem \ref{ribbon:contained:2:teo} extends \emph{as is} from $D^4$ to $H_g^4$ and furnishes a shadow $X$ of $G$ containing $S$. The shadow formula says that
$$\langle G \rangle = \sum_\sigma \langle X_\sigma \rangle.$$
The surface $S$ is the unique subsurface in $X$ with boundary equal to $L$: if there were another one $S'$, then $S+S'$ would be a non-trivial element in $H_2(X;\matZ_2)=\{e\}$. Therefore for every coloring $\sigma$ of $X$ extending that of $G$, the odd surface $S_\sigma$ coincides with $S$. Now Theorem \ref{order:teo} applies for each state $\langle X_\sigma \rangle$ and we are done.
\end{proof}

\subsection{Examples}
We show a couple of examples. The first one is pretty simple: let $S$ be a compact orientable surface with boundary. The four-manifold $S\times D^2$ is diffeomorphic to $H_g^4$ and its boundary is $\#_g(S^2\times S^1)$ with $g=1-\chi(S)$. The surface $S$ is a shadow for the link $L=\partial S$. It consists of one region $S$ with gleam zero.

If we color the components of $L$ with different colors, no coloring of $S$ can extend them: so there are no states, and $\langle L \rangle = 0$. In this case Theorem \ref{large:2:teo} gives no information. If we color each component of $L$ with the same $n$,
there is a single coloring $\sigma$ for $S$ extending it and we get
$$\langle L \rangle = \langle S_\sigma \rangle = \cerchio_n^{\chi(S)} = \big((-1)^n [n+1]\big)^{\chi(S)}.$$
When $n$ is odd, this function has a pole in $q=i$ of order $-\chi (S)$.
Therefore $S$ is the ribbon surface with smallest $-\chi(S)$ for $L$, and the lower bound given by Theorem \ref{large:2:teo} is sharp on these links. This proves Example \ref{easy:example}.

\begin{rem} \label{easy:was:rem}
It is in fact obvious that there cannot be any subsurface $S'\subset S\times D^2$ with $\partial S' = \partial S$ and $\chi(S') > \chi(S)$, because there is no map $S'\to S$ that sends $\partial S'$ homeomorphically to $\partial S$. (If we cap the surfaces we get a degree-one map from a lower-genus closed surface to a higher-genus one.)
\end{rem}

As another example we compute the Kauffman bracket of the framed knot $K\subset M_g$ drawn in Fig.~\ref{knot_ribbon:fig}, considered with its blackboard framing. We construct a shadow following the algorithm of Proposition \ref{Cost:prop}. To compute the gleam of the regions, we add some $\pm \frac 12$ around each crossing as follows:

\begin{center}
\includegraphics[width = 2.5 cm]{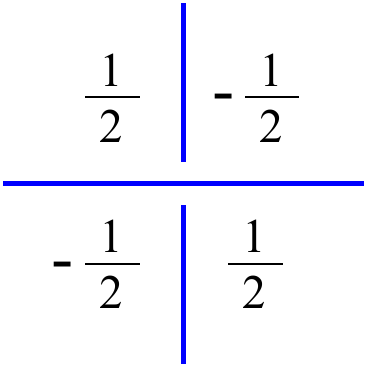}
\end{center}
and then we add all the contributions contained in each region, see \cite{CoThu, Tu}.

\begin{figure}
\begin{center}
\includegraphics[width = 11 cm]{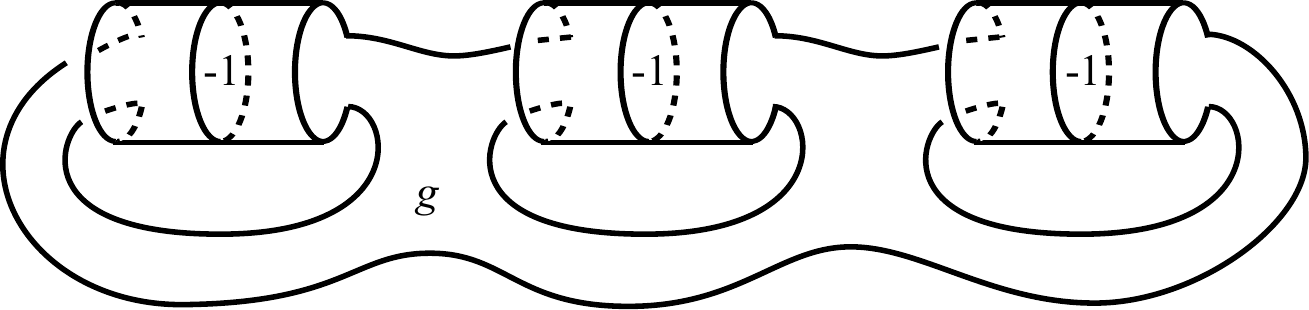}
\caption{A shadow for the knot $K$ drawn in Fig.~\ref{knot_ribbon:fig}. The $g$ discs have gleam $-1$ and the large region has gleam $g$.}
\label{ombra_conto:fig}
\end{center}
\end{figure}

As a result we get the shadow $X$ shown in Fig.~\ref{ombra_conto:fig}. The shadow $X$ has a large region $R$ with $\chi = -2g+1$ and gleam $g$, and $g$ discs $D_1,\ldots, D_g$ with gleam $-1$. 

We give $K$ the color 1. Every edge is circular and incident to $D_i, R, R$. If $D_i$ is colored by $c$, then $(c,1,1)$ must be an admissible triple: this holds only for $c=0,2$. Therefore each $D_i$ can be colored by wither $0$ or $2$. Thus a coloring $\sigma$ for $X$ is determined by a vector $\sigma = (\sigma_1,\ldots, \sigma_g) \in \{0,2\}^g$. 

The circular edges $e$ have $\chi(e)=0$ and hence do not contribute to the formula for $\langle X_\sigma \rangle$. Hence the only contributions come from the regions of $X$. Recall that a region $f$ with gleam $g$ and color $c$ contributes with a factor
$$\cerchio_c^{\chi(f)} q_f = \big((-1)^c[c+1]\big)^{\chi(f)} (\sqrt {-1})^{2gc}q^{-\frac{gc}2(c+2)}.$$
The large region $R$ contributes with
$$\cerchio_1^{\chi(R)} q_R = (-[2])^{-2g+1}(-1)^gq^{-\frac{3g}2} = (-1)^{1-g}\frac{q^{-\frac{3g}2}}{(q+q^{-1})^{2g-1}}.$$
A disc $D_i$ contributes according to its color $0$ or $2$ respectively as 
\begin{align*}
\cerchio_0 q_{D_i} & = 1, \\
\cerchio_2 q_{D_i} & = [3] q^4 = q^2+q^4+q^6.
\end{align*}
Set $|\sigma| = \sum_i\frac{\sigma_i}2$. We get:

\begin{align*}
\langle K \rangle & = \sum_\sigma \langle X_\sigma \rangle \\
& = \sum_\sigma (-1)^{1-g} \frac{q^{-\frac{3g}2}}{(q+q^{-1})^{2g-1}} \cdot (q^2+q^4+q^6)^{|\sigma|}
 \\
& = (-1)^{1-g} \frac{q^{-\frac{3g}2}}{(q+q^{-1})^{2g-1}} \sum_\sigma (q^2+q^4+q^6)^{|\sigma|} \\
& =  (-1)^{1-g} q^{-\frac{3g}2}\frac{(1+ q^2+q^4+q^6)^g}{(q+q^{-1})^{2g-1}}  \\
& = (-1)^{1-g}q^{\frac{3g}2}\frac{[4]^g}{[2]^{2g-1}}.
\end{align*}
Therefore:
$$\ord_i\langle K \rangle = g\cdot \ord_i[4] - (2g-1)\cdot \ord_i[2] = g-(2g-1) = 1-g.$$
This proves Example \ref{nodo:example}.

\section{The state-sum formula} \label{state-sum:section}
We prove here the shadow state-sum formula for $\langle G \rangle$, namely Theorems \ref{state-sum:teo} and \ref{state-sum2:teo}. Recall that $M_g = \#_g(S^2\times S^1)$ when $g\geqslant 1$, and we extend this notation by setting $M_0 = S^3$. We also set $H_0^4 = D^4$, so that $M_g = \partial H_g^4$ for all $g\geqslant 0$.

The shadow state-sum formula was first proved by Turaev \cite{Tu} in $S^3$ and hence extended by Costantino \cite{C:Jones} in $M_g$. We include here for completeness a proof that uses skein theory and avoids Reshetikhin-Turaev invariants.

\subsection{Fusions and sphere intersections}
\begin{figure}
\begin{center}
\includegraphics[width = 10 cm]{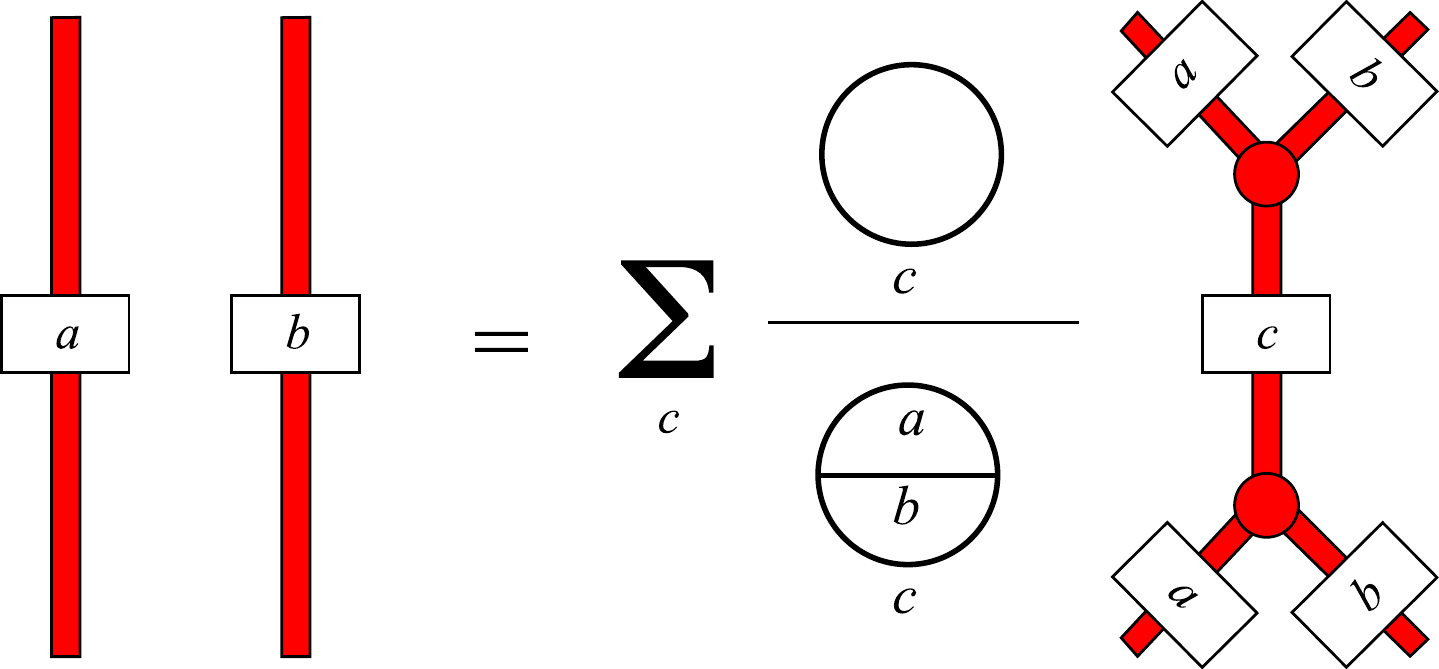}
\caption{The \emph{fusion rule}. Recall that all framings are orientable, \emph{i.e.}~they form an orientable surface that thickens the knotted trivalent graph. We suppose here that the two bands in the left are oriented coherently, so that the right knotted trivalent graph is also orientable.}
\label{fusion:fig}
\end{center}
\end{figure}

\begin{figure}
\begin{center}
\includegraphics[width = 12.5 cm]{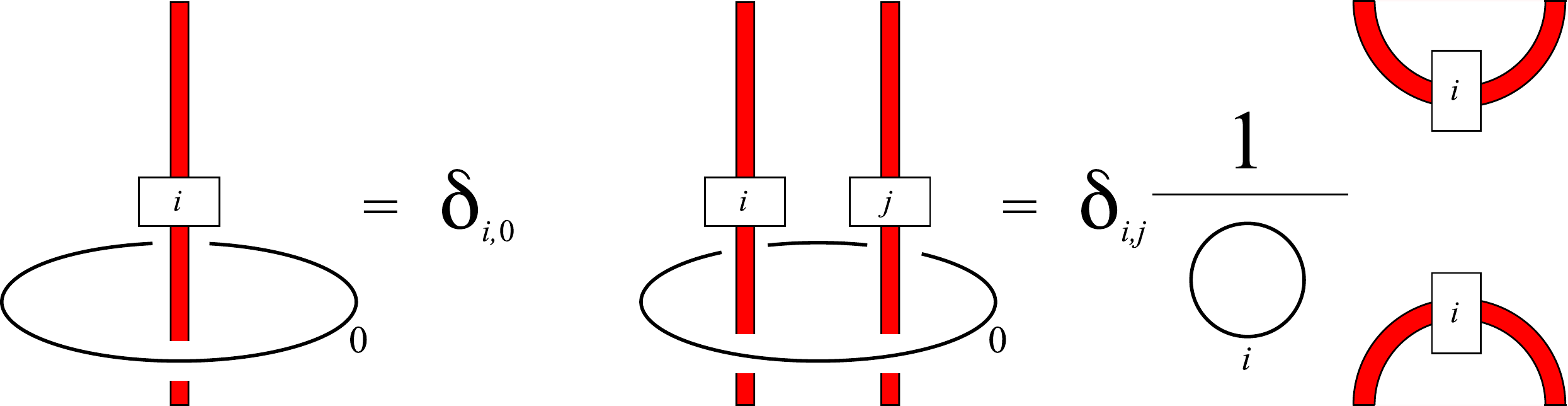}
\caption{Sphere intersection.}
\label{sphere:fig}
\end{center}
\end{figure}

We recall a couple of skein equalities.
The first is the well-known \emph{fusion rule} shown in Fig.~\ref{fusion:fig}, which takes place inside a ball, see \cite[Fig. 14.15]{L} 

A second kind of move is shown in Fig.~\ref{sphere:fig}-(left) and takes place in the neighborhood of a two-dimensional sphere $S$, drawn as a $0$-framed circle in the picture. If $G$ intersects $S$ transversely in exactly one point, then Fig.~\ref{sphere:fig}-(left) applies. The move says that if the edge of $G$ crossing $S$ has a positive coloring $i\geqslant 1$, then $G = 0$ as skeins. See \cite[Lemma 1]{YM} for a proof. Note that after applying the move we can surger along the sphere without affecting $\langle G \rangle$, see Remark \ref{canonical:rem}.

By combining the two moves we also get a third one that applies when $G$ intersects $S$ transversely into two points, see Fig.~\ref{sphere:fig}-(right). 

\subsection{Simple polyhedra that collapse onto graphs}
It might be non-obvious in general to determine whether a polyhedron collapses onto a graph. Luckily, on simple polyhedra there is a nice criterion.

\begin{prop}[Costantino] \label{collapses:prop}
Let $X$ be a connected simple polyhedron. The following facts are equivalent:
\begin{enumerate}
\item $X$ collapses onto a graph,
\item $X$ does not contain a simple polyhedron without boundary,
\item every coloring of $\partial X$ extends to finitely many colorings on $X$.
\end{enumerate}
\end{prop}
\begin{proof}
See \cite[Lemma 3.6]{C:Jones}.
\end{proof}

\begin{cor}
Let $X$ be a simple polyhedron that collapses to a graph. Every connected simple  subpolyhedron $X'\subset X$ also collapses onto a graph.
\end{cor}
\begin{proof}
The polyhedron $X$ does not contain any simple sub-polyhedron without boundary, hence $X'$ also does not.
\end{proof}

\begin{figure}
\begin{center}
\includegraphics[width = 9 cm]{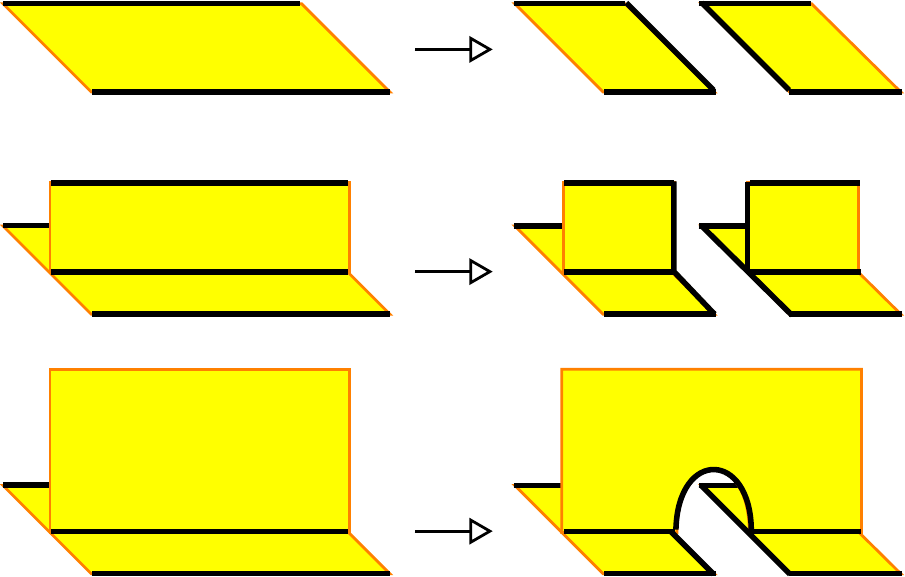}
\caption{A simple polyhedron $X$ that collapses onto a graph reduces to a finite union of atomic polyhedra after finitely many moves of this type. The bold exterior lines are portions of $\partial X$.}
\label{cut:fig}
\end{center}
\end{figure}

\begin{cor}
Each move in Fig.~\ref{cut:fig} transforms a simple polyhedron that collapses to a graph into one or two simple polyhedra that collapse to a graph.
\end{cor}

A simple polyhedron $X$ is \emph{atomic} if it is the cone over $\cerchio, \teta$, or $\tetra$, that is $X$ is as in Fig.~\ref{models:fig}-(3,2,1). We will use the following.

\begin{prop}
Let $X$ be a simple polyhedron that collapses onto a graph. The polyhedron reduces to a finite union of atomic polyhedra after a finite combination of moves as in Fig.~\ref{cut:fig}.
\end{prop}
\begin{proof}
We say that a region of $X$ is \emph{exterior} if it is incident to $\partial X$, and \emph{interior} otherwise. Suppose $X$ contains some interior regions. There is an edge $e$ that is adjacent to one interior region and to two exterior regions: if not, the interior regions would form a simple sub-polyhedron contradicting Proposition \ref{collapses:prop}. The move in Fig.~\ref{cut:fig}-(bottom) applied to $e$ transforms the interior region into an exterior one: after finitely many such moves we kill all the interior regions.

Now we can use Fig.~\ref{cut:fig}-(center) to cut every interior edge in two halves, and then Fig.~\ref{cut:fig}-(top) to cut every region into discs that are incident to $\partial X$ only in one arc or circle. We are left with atomic pieces.
\end{proof}

\subsection{Moves on shadows}
If we apply one of the moves of Fig.~\ref{cut:fig} to a shadow $X$ of some graph $G\subset M$, we get a new simple polyhedron $X'$ that can be interpreted as a shadow of some graph $G'$ in some manifold $M'$. We show this fact for each move.

We start by examining Fig.~\ref{cut:fig}-(top). The yellow strip thickens to a $D^3 \times [-1,1]$, with boundary $S^2\times [-1,1]$. The two-sphere $S = S^2\times 0$ intersects $G$ transversely into two points. Topologically, the move corresponds to surgerying $M$ along the two-sphere $S^2\times 0$ and modifying $G$ as in Fig.~\ref{sphere:fig}-(right). We get a new graph $G'$ inside a new manifold $M'$, with a new shadow $X'$. If $S$ is separating, these objects actually split into two components.

The move in Fig.~\ref{cut:fig}-(center) is analogous, the only difference being that now $S$ intersects $G$ in three points. The move in Fig.~\ref{cut:fig}-(bottom) is the fusion shown in Fig.~\ref{fusion:fig}.

\begin{rem} \label{cut:rem}
In the moves of Fig.~\ref{cut:fig}, some region $R\subset X$ is cut into two regions $R_1,R_2 \subset X'$. The gleams $g_1$ and $g_2$ of these new regions sum to give the gleam $g=g_1+g_2$ of $R$. The gleams of all the other regions of $X$ do not change.
\end{rem}

\subsection{The shadow formula}
We are now ready to prove the shadow formula. Recall that $M_g = \partial H_g^4$ and we use the convention $M_0 = S^3$ and $H_0^4 = D^4$.

\begin{teo}
Let $G$ be a colored framed knotted trivalent graph in $M_g$ and $X$ be a shadow for $G$, contained in $H_g^4$. We have
$$\langle G \rangle = \sum_\sigma \langle X_\sigma \rangle$$
where $\sigma$ varies among all colorings of $X$ extending that of $G$.
\end{teo}
\begin{proof}
We recall that
\begin{align} \label{state:eqn}
\langle X_\sigma \rangle = \frac{\prod_f \cerchio_f^{\chi(f)}q_f \prod_v \tetra_v^{\chi(v)} \prod_{v_\partial} \teta_{v_\partial}^{\chi(v_\partial)}}
{\prod_e \teta_e^{\chi(e)} \prod_{e_\partial} \cerchio_{e_\partial}^{\chi(e_\partial)}}.
\end{align}
The formula holds when $X$ is atomic with zero gleams: there is a single coloring $\sigma$ on $X$ extending that of $G$, and we get $\langle X_\sigma \rangle = \langle G \rangle$. To prove that, note that the contribution of every non-closed $e_\partial$ or $v_\partial$ cancels with the contribution of the incident $f$ or $e$. Therefore:
\begin{itemize}
\item if $G=\cerchio$ we get obviously $\cerchio$,
\item if $G=\teta$ everything cancels except $\teta_{v_\partial}^2 / \teta_e = \teta_e$,
\item if $G=\tetra$ everything cancels except $\tetra_v$.
\end{itemize}
Suppose now $X$ is atomic with arbitrary gleams. We modify the gleams using the following moves:
\begin{enumerate}
\item add a gleam $\pm 1$ on a region: this corresponds to a full twist of the corresponding framed edge of $G$;
\item add a gleam $\pm \frac 12$ to the three regions incident to an interior edge of $X$: this corresponds to a half-twist to each of the three edges of $G$ incident to a vertex of $G$.
\end{enumerate}
Using finitely many such moves we can reduce all gleams to zero. To show that,
color in green the regions having a half-integer (but non-integer) gleam. Recall that the framing of $G$ is orientable: this implies that every sub-circle $C\subset G$ intersects an even number of green faces, and it is easy to check that with moves (2) we can transform all gleams into integers. Then we reduce them to zero using (1).

Let $G'$ be obtained from $G$ by (1) or (2). We recall from \cite[Fig.~14.1 and 14.14]{L} that:
\begin{align*}
\langle G' \rangle & = (-1)^c q^{\mp\frac c2(c+2)} \langle G \rangle,\\
\langle G' \rangle & = (-1)^{\frac {a+b+c}2} q^{\mp\frac a4(a+2) \mp \frac b4(b+2) \mp \frac c4(c+2)} \langle G \rangle
\end{align*}
corresponding respectively to moves (1) and (2). In the formula (\ref{state:eqn}) the contribution of the phases 
$$q_f = (\sqrt{-1})^{2gc}q^{-\frac{gc}2(c+2)}$$
changes exactly in the same way: this proves the theorem for any atomic shadow $X$.

A more general $X$ decomposes into atoms via finitely many moves as in Fig.~\ref{cut:fig}. Let $n(X)$ be the number of moves necessary to atomize $X$: we prove the theorem by induction on $n(X)$. 

Pick a move transforming $X$ into a $X'$ with $n(X')<n(X)$. 
The polyhedron $X'$ is a shadow of some graph $G'$ in some manifold $M'$. The objects $X'$ and $M'$ may have two components, but the following arguments work anyway. We suppose by induction that the theorem holds for $X'$ and $G'$, and we prove it for $X$ and $G$.

Consider the move in Fig.~\ref{cut:fig}-(top). The pair $(M',G')$ is obtained from $(M,G)$ via the move shown in Fig.~\ref{sphere:fig}-(right), with $G'$ inheriting the coloring of $G$. Therefore
$$\langle G \rangle = \frac{1}{\cerchio_R} \langle G' \rangle$$
where $R$ is the yellow region that we have cut. There is an obvious correspondence between colorings of $X$ and $X'$, and the formula (\ref{state:eqn}) says that for each coloring $\sigma$ we have
$$\langle X_\sigma \rangle  = \frac 1{\cerchio_R} \langle X'_\sigma \rangle.$$
(We use here Remark \ref{cut:rem} to show that the phases of $X_\sigma$ and $X'_{\sigma}$ are the same.)
The theorem holds for the pair $(X', G')$, and hence holds also for $(X,G)$.

The move in Fig.~\ref{cut:fig}-(center) is treated analogously. Using a fusion and Fig.~\ref{sphere:fig} we find easily that
$$\langle G \rangle = \frac{1}{\teta_e} \langle G' \rangle $$
where $e$ is the edge cut in Fig.~\ref{cut:fig}-(center). There is an obvious correspondence between colorings of $X$ and $X'$, and for  each such coloring $\sigma$ we have
$$\langle X_\sigma \rangle  = \frac 1{\teta_e} \langle X'_\sigma \rangle.$$

Finally, the move in Fig.~\ref{cut:fig}-(bottom) is a fusion. The fusion formula says
$$\langle G \rangle = \sum_c \frac{\cerchio_c}{\teta_{a,b,c}} \langle G'_c \rangle$$
where the coloring $G'_c$ on $G'$ varies on the new edge $c$. Every coloring of $X$ induces one of $X'$ and we get
$$\langle X_\sigma \rangle  = \frac {\cerchio_c}{\teta_{a,b,c}} \langle X'_\sigma \rangle.$$
This proves the theorem.
\end{proof}

\section{Open questions} \label{questions:section}
A list of stimulating open questions is contained in the last section of Eisermann's paper \cite{Eis}, which is overall very nice and enjoyable to read. Here we add more questions to that list.

\subsection{Ribbon genus of knots in $S^3$}
We have seen that the Kauffman bracket $\langle K \rangle$ of a knot $K\subset \#_g(S^2\times S^1)$ may produce non-trivial (sometimes sharp) lower bounds for the ribbon genus of $K$, see Examples \ref{nodo:example} and \ref{easy:example}.

The situation in $S^3$ is more disappointing because of the following:

\begin{prop}
Let $L\subset S^3$ be a colored framed link. If at least one component of $L$ has an odd coloring, the bracket $\langle L \rangle$ vanishes at $q=i$.
\end{prop}
\begin{proof}
We prove it by induction on the maximum color $c$ on $L$. If $c=1$, this is the standard case: we choose a diagram for $L$ and use the first Kauffman bracket relation to transform $\langle L \rangle$ into a linear combination of unlinks with coefficients in $\matZ[A^{\pm 1}]$. The second bracket relation says that the bracket of each unlink vanishes at $q=i$.

\begin{figure}
\begin{center}
\includegraphics[width = 6 cm]{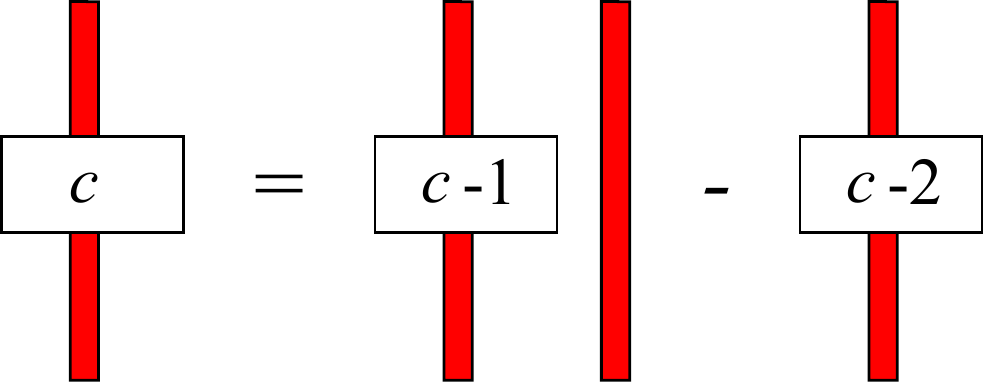}
\caption{A skein move on a framed colored knot.}
\label{Cheb:fig}
\end{center}
\end{figure}

If some component $K$ of $L$ has a color $c>1$, we modify $K$ via the well-known skein move shown in Fig.~\ref{Cheb:fig} that takes place in a solid torus neighborhood of $K$ and is an immediate consequence of Fig.~\ref{JW:fig}. Each of the new two addenda is a colored link with at least one odd-colored component. We perform this move on all components with maximum color $c$ and we conclude by induction.
\end{proof}

Therefore $\ord_i\langle K \rangle \geqslant 1$ for every odd-colored knot $K\subset S^3$ and if we apply Theorem \ref{large:teo} to $\langle K \rangle$ we get no relevant information. One could try however to choose a colored knotted trivalent graph $G\subset S^3$ containing the knot $K$ as its odd sublink. We do not know if some relevant information may be obtained for $K$ in that case:

\begin{quest}
Is there a colored framed knotted trivalent graph $G\subset S^3$ whose odd sub-link $K\subset G$ is a knot, such that
$$\ord_i\langle G \rangle +\frac r2 \leqslant 0\ ?$$
\end{quest}
One such example would imply that $K$ is not ribbon. 
\begin{rem}
As far as we know, it might be that $\ord_i\langle G \rangle + \frac r2 >0$ for all colored trivalent $G\subset S^3$ having a non-empty odd sub-link. See for instance \cite{C:integral} where it is shown that $\langle G \rangle$ is a polynomial up to a little renormalization.
\end{rem}

More generally, we do not know if by passing from links to graphs we gain more obstructions for the existence of ribbon surfaces, because we tested only very few examples. Computing the Kauffman bracket of a colored knotted trivalent graph $G\subset S^3$ by hand can be tedious: it would be nice to have a computer program where the user can draw a diagram of $G$ and get $\ord_i \langle G \rangle$ as a result. We have computed by hand a couple of examples (the Hopf link and the trefoil knot with an additional arc) and found no improvement there.

\subsection{More manifolds}
The notion of ribbon surface applies to any kind of 3-manifold $M$, but the Jones polynomial does not. To define $\langle K \rangle$ as a rational function we need the Kauffman space $K(M)$ to be one-dimensional.

\begin{quest}
For which closed 3-manifolds $M$ the space $K(M)$ is one-dimensional?
\end{quest}

When $K(M)$ is not one-dimensional, quantum invariants survive only at the roots of unity: these are the well-known Reshetikhin-Turaev-Witten invariants. These invariants can also be calculated using shadows, so it might be that some of the techniques used here extend to that context:

\begin{quest}
Can we relate the ribbon genus of a link to RTW invariants, for instance by taking roots of unity $q$ converging to $q \to i$? Does the fact that a knot is ribbon influence the asymptotic of the RTW invariants as $q \to i$?
\end{quest}

\subsection{Ribbon surfaces and shadows}
We have discovered that being contained in a shadow is a property that lies in the middle, between being ribbon and being homotopically ribbon. It is then natural to ask  Question \ref{intro:quest}, which splits into two questions. Let $S$ be a properly embedded surface in $D^4$:

\begin{quest}
If $S$ is contained in a shadow, is it ribbon? 
\end{quest}

\begin{quest}
If $S$ is homotopically ribbon, is it contained in a shadow?
\end{quest}

\subsection{Slice-ribbon conjecture} \label{slice-ribbon:subsection}
The famous slice-ribbon conjecture states that a knot in $S^3$ is \emph{slice} (\emph{i.e.}~it bounds a smooth disc in $D^4$) if and only if it is ribbon. It is worth mentioning that this conjecture extends naturally at least in three ways: from knots to links, from discs to more general surfaces, and also from $S^3$ to more general 3-manifolds.
Since we have not seen it in the literature, we state this three-fold generalization as a question:

\begin{quest}
Let $M$ be any 3-manifold.
Let $L\subset M$ be a link in $M=M\times 0$ that bounds a compact properly embedded surface $S\subset M\times [0,1]$. Does $L$ bound a \emph{ribbon} surface $S'$ diffeomorphic to $S$?
\end{quest}

We may define the \emph{slice genus} $g_s(K)$ of a null-homologous knot $K\subset M$ as the smallest genus of a properly embedded orientable surface $S$ in $M\times [0,1]$ with $\partial S = K$, and the \emph{ribbon genus} $g_r(K)$ as the smallest genus of a ribbon surface $S$ for $K$. When $M=S^3$ these are the standard slice and ribbon genera, since every surface in $D^4$ can be pushed inside $S^3 \times [0,1]$.

Of course we have $g_r(K) \geqslant g_s(K)$, and the previous question specializes to the following.

\begin{quest}
Does the equality $g_r(K) = g_s(K)$ hold for every possible null-homologous knot $K\subset M$ in every 3-manifold $M$?
\end{quest}

The lower bounds for the ribbon genus proved in this paper might in principle be used to find a counterexample in $M=\#_g(S^2\times S^1)$.

\subsection{Other roots of unity}
Let $X$ be a simple spine of a 3-manifold $M$. Roughly, a simple spine is just a shadow with all gleams zero: spines are used for instance in Turaev-Viro invariants \cite{TV}. For instance, $X$ might be the dual of an ideal triangulation for $M$.

A coloring $\sigma$ for $X$ gives rise to a rational function $\langle X_\sigma \rangle$ that may have poles in $q=0,\infty$, and at some roots of unity. The coloring $\sigma$ defines a \emph{spinal surface} $F_\sigma \subset M$, and we get 
$$\ord_0 \langle X_\sigma \rangle\geqslant -\chi(F_\sigma)$$ 
by a nice result of Frohman and Kania-Bartoszynska \cite{FK:content} that connects quantum invariants near $q=0$ to normal surfaces theory. This result was used extensively for instance in \cite{CM}. Here we have proved that
$$\ord_i \langle X_\sigma \rangle \geqslant \chi (S_\sigma)$$
where $S_\sigma$ is the odd surface contained in $X$. Note that the two inequalities concern different surfaces, and have opposite signs!

\begin{quest}
Do we get any similar inequalities for $\ord_q\langle X_\sigma \rangle$ when $q$ is a root of unity different from $\pm i$?
\end{quest}

\end{document}